\documentclass[10pt]{article}
 \usepackage[T1]{fontenc}
 \usepackage{latexsym,amssymb,amsmath}
\usepackage{pstricks,lscape,rotating}
\usepackage{fancyhdr}

 \newcommand{\G}{\widetilde{G}}

 \newcommand{\Galg}{\mathbf{G}}
 
 \newcommand{\Balg}{\mathbf{B}}
 
 \newcommand{\Lalg}{\mathbf{L}}

 \newcommand{\Ualg}{\mathbf{U}}
 \newcommand{\Talg}{\mathbf{T}}

 \newcommand{\Nn}{\operatorname{N}}
 \newcommand{\Stab}{\operatorname{Stab}}

 \newcommand{\Tr}{\operatorname{Tr}}
 \newcommand{\Gf}{{G}}

 \newcommand{\s}{\sigma}

 \newcommand{\Res}{\operatorname{Res}}

 \newcommand{\Cen}{\operatorname{C}}
 \newcommand{\Cl}{\operatorname{Cl}}
 
 \newcommand{\St}{\operatorname{St}}
 \newcommand{\rank}{\operatorname{rk}}

 \newcommand{\semi}[1]{\rtimes\langle\,#1\,\rangle}
 
 \newcommand{\cyc}[1]{\langle\,#1\,\rangle}

 \newcommand{\F}{\mathbb{F}}

 \newcommand{\Q}{\mathbb{Q}}
 \newcommand{\Z}{\mathbb{Z}}

 \newcommand{\wG}{\widetilde{\Galg}}

 \newcommand{\Irr}{\operatorname{Irr}}

\title{\bf On Lusztig's conjecture for connected and disconnected exceptional groups}

\author{Olivier \textsc{Brunat}\\
 brunat@math.rwth-aachen.de}

 \newtheorem{lemme}{Lemma}[section]
 
 \newtheorem{theoreme}{Theorem}[section]
 
 \newtheorem{proposition}{Proposition}[section]
 \newtheorem{corollaire}{Corollary}[section]
 \newtheorem{remarque}{Remark}[section]

 \newtheorem{conjecture}{Conjecture}[section]
 \newenvironment{preuve}[1][]{\noindent {\bf{Proof }}{---\ }}{\hfill
 \nopagebreak $\Box$
 \newline
}

\newenvironment{changemargin}[2]{\begin{list}{}{%
\setlength{\topsep}{0pt}%
\setlength{\leftmargin}{0pt}%
\setlength{\rightmargin}{0pt}%
\setlength{\listparindent}{\parindent}%
\setlength{\itemindent}{\parindent}%
\setlength{\parsep}{0pt plus 1pt}%
\addtolength{\leftmargin}{#1}%
\addtolength{\rightmargin}{#2}%
}\item }{\end{list}}

\begin{document}
\maketitle

\begin{changemargin}{2cm}{2cm}
\begin{center}{\sc abstract}
\end{center}
Lusztig conjectured that the almost characters of a finite reductive
group are up to a scalar the same as the characteristic functions of
the rational character sheaves defined on the corresponding algebraic
group. We propose in this paper to verify this conjecture for the
uniform almost characters of the Suzuki and the Ree groups and to show
that a conjecture of the same type holds for the disconnected case of
type~$B_2$. We then state a conjecture on the possible values of the
uniform almost characters on unipotent elements in the disconnected
case of type~$F_4$. 
\end{changemargin}

\section{Introduction}

Let~$\Galg_0$ be a connected reductive group defined over the finite
field with~$q$ elements and~$F_0:\Galg_0\rightarrow\Galg_0$ such 
that~$F_0^{\delta}$ is a Frobenius map for some positive
integer~$\delta$. We denote by~$\Galg_0^{F_0}$ and~$\Galg_0^{F_0^{\delta}}$ the
finite groups of fixed points under~$F_0$ and~$F_0^{\delta}$ respectively.
The group~$\Galg_0^{F_0^{\delta}}$ is~$F_0$-stable and the restriction of~$F_0$
to~$\Galg_0^{F_0^{\delta}}$ (denoted by the same symbol) is an automorphism
of~$\Galg_0^{F_0^{\delta}}$. We are interested in the irreducible
characters of the extension~$\G_0=\Galg_0^{F_0^{\delta}}\semi{F_0}$. A
relevant example in which this situation occurs is when~$\Galg_0$
is a simple algebraic group of type~$B_2$,~$G_2$ or~$F_4$ and~$F_0$ is
the isogeny that defines the Suzuki and the Ree groups; in these
cases,~$F_0^2$ is a Frobenius map and~$\Galg_0^{F_0^2}$ is a simple group 
(except for the first value of the parameter).
In~\cite{brunat} and~\cite{brunat1}, we compute the character table
of~$\G_0$ in the cases of type~$B_2$ and~$G_2$ respectively. But the employed methods
are essentially \emph{ad hoc} and cannot be used for type~$F_4$. 

When the algebraic group is connected, Lusztig has shown that the
computation of the irreducible characters of~$\Galg_0^{F_0}$ is reduced to the
determination of the \emph{almost characters} of~$\Galg_0^{F_0}$ (this is a
family of class
functions on~$\Galg_0^{F_0}$ defined using the generalized characters of
Deligne-Lusztig)~\cite{Lust}. To obtain the values of these functions, he develops
the theory of character sheaves on~$\Galg_0$: to every~$F_0$-stable
character sheaf on~$\Galg_0$ we can associate in a natural way a class
function on~$\Galg_0^{F_0}$, the so-called characteristic function of the
sheaf~\cite{Lustsheave}. He conjectures that the characteristic functions of character
sheaves coincide (up to a scalar) with the almost characters~\cite[V]{Lustsheave}.
With some assumptions on the characteristic and on the center
of~$\Galg_0$, Shoji showed that this conjecture holds when~$F_0$ is a
Frobenius map~\cite{shoji}. Recently, Bonnaf\'e proved that Lusztig's
conjecture holds in type~$A_n$ when the center of~$\Galg_0$ is not
connected~\cite{bonnafe}. 

The advantage of this approach is that the characteristic functions of
$F_0$-stable character sheaves can be computed using a uniform algorithm
described by Lusztig in~\cite[V.\S24]{Lustsheave}; the input of this
algorithm is the generalized Springer correspondence
introduced by Lusztig~\cite{Lustinter}.

Recently, Lusztig started a series of papers in which he
develops a theory of character sheaves for disconnected groups
and defines a generalized Springer correspondence in this
case~\cite{Lustdisconnected}; we do not know at present whether Lusztig's
algorithm is valid in this situation.

In this note, we propose to apply this theory to the groups~$\G_0$
interpreted as a fixed-point subgroup of a disconnected reductive
group~$\widetilde{\Galg}=\Galg\semi\tau$ (here~$\Galg$ is a connected
reductive group constructed from~$\Galg_0$ and~$\tau$ is an automorphism of 
finite order of~$\Galg$) under a generalized Frobenius~$F$, using a method 
developed
by Digne in~\cite{Digne}. In the case that~$\Galg_0$ is of type~$B_2$
or~$F_4$, the element~$\tau$ (viewed as an element of~$\widetilde{\Galg}$)
is unipotent and we thus can determine the generalized Springer
correspondence on the coset~$\Galg\tau$. We can then apply Lusztig's
algorithm \emph{a priori} as in the connected case; we show that Lusztig's
conjecture holds for the disconnected type~$B_2$ (using the explicit
table in~\cite{brunat}) and we give the possible values of the uniform almost
characters on unipotent elements for the type~$F_4$.

The paper is organized as follows: in Section~\ref{sec2} we determine
the classes of~$\widetilde{\Galg}$ from those of~$\Galg_0$.
In Section~\ref{sec3},
we apply Lusztig's algorithm to~$\Galg_0$ (in the case of
types~$B_2$,~$G_2$ and~$F_4$) with the twisted map~$F_0$ that defines
the Suzuki and Ree groups and verify that Lusztig's conjecture holds in those
cases for the uniform almost characters on the unipotents elements. 
In Section~\ref{sec4}, we verify Lusztig's conjecture in the
disconnected case of type~$B_2$ and give a conjecture for the case~$F_4$.

\section{The group~$\G_0$}\label{sec2}

In this section,~$\Galg_0$ is a connected reductive group defined
over~$\F_q$ (where~$q$ is a power of~$2$) and~$F_0$ is a map
on~$\Galg_0$ such that~$F_0^2$ is a Frobenius map.

\subsection{Restriction of scalars}\label{rescal}

Now we will interpret the group~$\G_0=\Galg_0^{F_0^2}\semi{F_0}$
as a fixed-point subgroup under an endomorphism~$F$ of an algebraic 
group~$\widetilde{\Galg}$;
we use the construction developed in~\cite{Digne}. We set~$\Galg=\Galg_0\times\Galg_0$
and~$\widetilde{\Galg}=\Galg\rtimes\cyc{\tau}$, where~$\tau$
is the algebraic automorphism of~$\Galg$ defined
by~$(g_1,g_2)^{\tau}=(g_2,g_1)$. The group~$\widetilde{\Galg}$
is a non-connected reductive group with identity
component~$\Galg$. Let~$F:\wG\rightarrow\wG$ be the homomorphism defined
by~$F(g_1,g_2)=(F_0(g_2),F_0(g_1))$ and~$F(\tau)=\tau$.
Let~$\Talg_0$ be a maximal rational torus of~$\Galg_0$ contained in a 
rational Borel subgroup of~$\Galg_0$.
The group~$\Balg=\Balg_0\times\Balg_0$ is an $F$-rational
Borel subgroup of~$\Galg$ containing~$\Talg= \Talg_0\times\Talg_0$
as maximal $F$-rational torus. The groups~$\Talg$ and~$\Balg$ are
$\tau$-stable, thus~$\tau$ is~\emph{quasi semisimple}. It is not
difficult to show that~$\tau$ is in fact \emph{quasi-central}
\cite[1.15]{DMnonconnexe}. We put~$\widetilde{\Balg}=\Balg\semi{\tau}$
and~$\widetilde{\Talg}=\Talg\semi{\tau}$; we
have~$\widetilde{\Balg}=\Nn_{\wG}(\Balg)$
and~$\widetilde{\Talg}=\Nn_{\wG}(\Talg,\Balg)$. Thus~$\widetilde{\Talg}$
and~$\widetilde{\Balg}$ are respectively rational "torus" and rational "Borel"
in the sense of Digne and Michel~\cite[1.2]{DMnonconnexe}.

\begin{proposition}\label{isogroupe} We have the following isomorphisms:
$$\begin{array}{lcl}
{\Galg}^{\tau}&\simeq&\Galg_0,\\
\left({\Galg}^{\tau}\right)^{F}&\simeq&\Galg_0^{F_0},\\
{\Galg}^F&\simeq&\Galg_0^{F_0^2},\\
{\widetilde{\Galg}}^F&\simeq&\Galg_0^{F_0^2}\semi{F_0}.
\end{array}$$
\end{proposition}

\begin{preuve}{}
The isomorphisms are given
by~$\varphi:\Galg_0\rightarrow{\Galg}^{\tau},\ g\mapsto (g,g)$
and~$\varphi':\Galg_0^{F_0^2}\rightarrow{\Galg}^F,\ g\mapsto
(g,F_0(g))$. Moreover we have~$\varphi'\circ F_0=\tau\circ\varphi'$,
it follows that~$F|_{{\Galg}^F}$ acts on~${\Galg}^F$
as~$F_0|_{\Galg_0^{F_0^2}}$ on~$\Galg_0^{F_0^2}$ and the result holds.

\end{preuve}

\subsection{Conjugacy classes}\label{secconjclass}

We will now describe the conjugacy classes of the groups appearing
in~\S\ref{rescal}. We are interested in the conjugacy classes
of elements belonging to the coset~$\Galg\tau$.\medskip

\begin{lemme}\label{critereconj}
Two elements~$(g_1,g_2)\tau$
and~$(g_1',g_2')\tau$ of~$\Galg\tau$ are conjugate in~$\wG$ 
if and only if~$g_1g_2$ and~$g_1'g_2'$ are conjugate in~$\Galg_0$.
\end{lemme}

\begin{preuve}{}First recall that~$(g_1,g_2)\tau$
and~$(g_1',g_2')\tau$ are conjugate in~$\wG$ if and only if
they are conjugate by an element in~$\Galg$. Then there exists
$h=(h_1,h_2)\in\Galg$ such that~$h_1g_1=g_1'h_2$ and~$ h_2g_2=g_2'h_1$.
It follows that~$g_1g_2$ and~$g_1'g_2'$ are conjugate. Conversely
if~$g_1g_2=hg_1'g_2'h^{-1}$ for some~$h\in\Galg_0$, we have:
$$\begin{array}{lcl}
(g_1,1)(g_2,g_1)\tau(g_1,1)^{-1}&=&(g_1g_2,1)\tau\\ 
&=&(hg_1'g_2'h^{-1},1)\tau\\
&=&(hg_1',h)(g_2',g_1')\tau(hg_1',h)^{-1}\\
\end{array}$$
Since~$\tau(g_1,g_2)\tau\tau^{-1}=(g_2,g_1)\tau$ the result follows.

\end{preuve}
In the following, if~$H$ is a subgroup of~$\Galg_0$, we denote by~$\mu(H)=
\{(h,h)\in \Galg\ |\ h\in H\}\subseteq\Galg$.\medskip

\begin{proposition}\label{bijection}
The map~$f:\Galg_0\mapsto\Galg\tau,\ x\mapsto(x,1)\tau$
induces a bijection from the conjugacy classes of~$\Galg_0$
to the classes in the coset~$\Galg\tau$ and we have
$$\Cen_{\wG}(f(x))=\mu(\Cen_{\Galg_0}(x)).\cyc{(x,1)\tau}.$$
Moreover~$f$ induces a bijection between the unipotent classes
of~$\Galg_0$ and the unipotent classes in the coset~$\Galg\tau$.
\end{proposition}

\begin{preuve}{}
We first remark that every conjugacy class in the coset~$\Galg\tau$
has a representative of the form~$f(g)$ for some~$g\in\Galg_0$
because~$(x_1,x_2)\tau$ is conjugate to~$f(x_1x_2)$. The bijection
between the classes is now a consequence of Lemma~\ref{critereconj}.
Moreover we have:
$$\left((x_1,x_2)\tau\right)^{2m}=\left((x_1x_2)^m,(x_2x_1)^m\right).$$
We denote by~$d$ the order of~$(x_1,x_2)\tau$ and by~$d'$ the one of~$x_1x_2$. 
Since~$x_1x_2$ and~$x_2x_1$ are conjugate, they have the same order. Thus we have~$((x_1,x_2)\tau)^{2d'}=1$ 
and it follows that~$d$ divides~$2d'$. Moreover~$2$ divides~$d$, 
then there exists a positive integer~$k$ such that~$d=2k$. Thus~$k$ divides~$d'$. Conversely we have~$1=((x_1,x_2)\tau)^{2k}$; 
it follows that~$(x_1x_2)^k=1$ and~$d'$ divides~$k$. We have~$k=d'$ and we have to prove that
$$|\cyc{(
(x_1,x_2),\tau)}|=2|\cyc{x_1x_2}|.$$ Recall that the characteristic of
the considered groups is~$2$, it follows that the bijection induced
by~$f$ is compatible with the unipotent classes. 

Let~$x\in\Galg_0$, we then
have~$\mu(\Cen_{\Galg_0}(x)).\cyc{(x,1)\tau}\subseteq\Cen_{\Galg}(f(x))$. 
Conversely if~$(g_1,g_2)(x,1)\tau=(x,1)\tau(g_1,g_2)$ we
have~$(g_1x,g_2)\tau=(xg_2,g_1)\tau$ and~$g_1\in\Cen_{\Galg_0}(x)$.
If~$(g_1,g_2)\tau(x,1)\tau=(x,1)\tau(g_1,g_2)\tau$, then $(g_1,g_2 x)=(x
g_2,g_1)$ and~$g_2\in\Cen_{\Galg_0}(x)$ as required.

\end{preuve}

We now discuss the $F$-rational classes of~$\wG$ in the coset~$\Galg\tau$.
We first recall a relevant result which is a direct consequence of
Lang's theorem:\medskip

\begin{theoreme}\label{gk} Let~$\Galg$ be a connected reductive
group defined over~$\F_q$ and~$F$ a generalized Frobenius
map on~$\Galg$. Suppose that~$\Galg$ acts transitively on a
set~$X$, suppose that there exists~$F':X\rightarrow X$ such
that~$F'(gx)=F(g)F'(x)$ for every~$g\in\Galg$ and~$x\in X$
and suppose finally that~$\Stab_{\Galg}(x)$ is a closed
subgroup of~$\Galg$. Then~$X^{F'}$ is non-empty and if~$x_0\in X^{F'}$, 
we set
$$A(x_0)=\Stab_{\Galg}(x_0)/\Stab_{\Galg}(x_0)^{\circ}$$ and by~$\pi$
the canonical projection of~$\Stab_{\Galg}(x_0)$ onto~$A(x_0)$. Then the
map~$F$ naturally induces an automorphism on~$A(x_0)$, denoted by the
same symbol, and the~$\Galg^F$-orbits of~$X^{F'}$ are in bijection
with the~$F$-classes of~$A(x_0)$ (i.e., the orbits on~$A(x_0)$ for
the operation~$g.h=g^{-1}hF(g)$). Moreover the correspondence is such
that if~$x=gx_0$ then we associate to the~$\Galg^F$-class of~$x$
the~$F$-class of~$\pi(g^{-1}F(g))$.
\end{theoreme}

Let~$x\tau\in\wG$, we denote
by~$A(x\tau)=\Cen_{\Galg}(x\tau)/{\Cen_{\Galg}(x\tau)}^{\circ}$.
If~$x\tau$ is rational, we denote by~$F$ the automorphism
of~$A(x\tau)$ induced by~$F$. We have:

\begin{proposition}\label{classesext} Every rational class
on~$\Galg\tau$ has an~$F$-stable representative. Moreover the map~$f$
defined in Proposition~\ref{bijection} gives a bijection between
the~$F_0$-rational classes of~$\Galg_0$ and the~$F$-classes in the
coset~$\Galg\tau$. Finally, the conjugacy classes of~$\Galg_0^{F_0}$
are in bijection with the conjugacy classes of~$\wG^{F}$ in the
coset~$\Galg^F\tau$. Moreover if~$g$ is a representative of a class
of~$\Galg_0^{F_0}$ and if~$\tilde{g}$ is a representative of the corresponding
class on~$\Galg^F\tau$ then we have
$$\begin{array}{lcl}
\Cen_{\widetilde{\Galg}^F}(\widetilde{g})&\simeq&\mu(\Cen_{\Galg_0^F}(g)).\cyc{(g,1)\tau},\\
|\Cen_{\widetilde{\Galg}^F}(\widetilde{g})|&=&2|\Cen_{\Galg_0^F}(g))|.
\end{array}
$$

\end{proposition}

\begin{preuve}{} 
The map~$f$ gives a bijection between
the~$F_0$-rational classes of~$\Galg_0$ and the~$F$-classes in the
coset~$\Galg\tau$ because~$F\left((g,1)\tau\right)$ is
conjugate to~$(g,1)\tau$ if and only if~$F_0(g)$ and~$g$ are conjugate.
Let~$g_0\in\Galg_0^{F_0}$, then we have
$$F((g_0,1)\tau)=(g_0,1)^{-1}(g_0,1)\tau(g_0,1).$$ 
Since~$(g_0,1)\in\Galg$, we can find~$u\in\Galg$ such
that~$(g_0,1)=u^{-1}F(u)$ (we apply Lang's property in the
connected group~$\Galg$ with the generalized Frobenius~$F$). The
element~$u(g_0,1)\tau u^{-1}$ is rational and we
have $$u\Cen_{\Galg}((g_0,1)\tau)u^{-1}=\Cen_{\Galg}(u(g_0,1)\tau
u^{-1}).$$ Let~$g,h\in\Cen_{\Galg}((g_0,1)\tau)$, we then have
$$uh^{-1}u^{-1}ugu^{-1}F(u)F(h)F(u^{-1})=uh^{-1}g(g_0,1)F(h)(g_0,1)^{-1}u^{-1}.$$
Consequently the $F$-classes of~$\Cen_{\Galg}(u(g_0,1)\tau
u^{-1})$ are in bijection with the $^{(g_0,1)}F$-classes
of~$\Cen_{\Galg}((g_0,1)\tau)$. We can apply 
Theorem~\ref{gk} with~$X=\Cl((g_0,1)\tau)$ and~$F'$ being
the restriction of~$F$ to this set; we take~$u(g_0,1)\tau u^{-1}$ as
rational element. It follows that the~$\wG^F$-classes of~$\Cl(f(x_0))^F$
are in bijection with the~$^{(g_0,1)}F$-classes
of~$A(f(x_0))$. Using Proposition~\ref{bijection}, we
have~$\Cen_{\Galg}(f(g_0))=\mu\left(\Cen_{\Galg_0}(g_0)\right)$.
For every~$(h,h)\in\Cen_{\Galg}(f(g_0))$, we then have~$$^{(g_0,1)}F(h,h)=(^{g_0}F_0(h),F_0(h)),$$
but since~$F_0$ acts on~$\Cen_{\Galg_0}(g_0)$ we have~$F_0(h)\in\Cen_{\Galg_0}(g_0)$ and it follows that
the operation of~$^{(g_0,1)}F$ on~$\Cen_{\Galg}(f(g_0))$ is the same as the one of~$F$ on~$\Cen_{\Galg}(f(g_0))$.
Since~$\mu$ is an isomorphism of algebraic
groups, it follows that~$A(f(g_0))$ is isomorphic
to~$A_{\Galg_0}(g_0)=\Cen_{\Galg_0}(g_0)/\Cen_{\Galg_0}(g_0)^{\circ}$ and we have 
$$F\circ\mu=\mu\circ F_0,$$ thus
the~$F_0$-classes of~$A_{\Galg_0}(g_0)$ are in bijection with
the~$F$-classes of~$A(f(g_0))$. This proves that the
classes of~$\Galg_0^F$ are in correspondence with the classes
of~$\wG^F$ in the coset~$\Galg\tau$. We set~$\pi$ (resp.~$\pi_0$) to be
the canonical projection of~$\Cen_{\Galg}(f(g_0))$ onto~$A(f(g_0))$
(resp.~$\Cen_{\Galg_0}(g_0)$ onto~$A_{\Galg_0}(g_0)$). Let~$\pi_0(a)$
be a representative of a~$F_0$-class of~$A_{\Galg_0}(g_0)$ and
let~$\widetilde{a}$ be a corresponding class representative
of~$\Galg_0^{F_0}$ given in Theorem~\ref{gk}. We then have
$$\Cen_{\Galg_0^{F_0}}(\widetilde{a})\simeq\{h\in\Cen_{\Galg_0}(g_0)\ |\
F(h)=a^{-1}ha\}.$$ Moreover the element~$\pi(a,a)=(\pi_0(a),\pi_0(a))$
is a representative of the~$^{(g_0,1)}F$-class
of~$A(f(g_0))$. Thus
$$\Cen_{\Galg^{F}}({\widetilde{a,a}})\simeq\{h\in\Cen_{\Galg}(f(g_0))\
|\ F(h)=a^{-1}ha\},$$ where~$(\widetilde{a,a})$
is a representative of the class of~$\wG^F$ corresponding
to~$\pi(a,a)$. 
Hence~$\Cen_{\Galg^F}(\widetilde{a,a})$ is isomorphic to
$\mu(\Cen_{\Galg_0^F}(\widetilde{a}))$ and the result is proven.

\end{preuve}

\begin{remarque}\label{shin}
The element~$u$ in the preceding proof can be explicitly
chosen as follows: since~$g_0$ belongs to the connected
group~$\Galg_0$, using Lang's property in this group, we
can find~$x\in\Galg_0$ such that~$g_0=x^{-1}F_0^2(x)$. We then
have $$(g_0,1)=(x,F_0(x))^{-1}F(x,F_0(x)).$$ It follows that
$y=\left(F_0(x)x^{-1},F_0(F_0(x)x^{-1})\right)\tau$ is a rational
element in the rational class of~$(g_0,1)\tau$. In~\cite{DM}
Digne and Michel define the Shintani correspondence: this is a
bijection~$N_{F/F^2}$ between the conjugacy classes of~$\Galg_0^{F_0}$
and the conjugacy classes of~$\Galg_0^{F_0^2}\semi{F_0}$
in the coset~$\Galg_0^{F_0^2}F_0$. We then remark that
$$y=\varphi'\left((N_{F/F^2}(g_0))^{-1}\right),$$ where~$\varphi'$ is
the explicit isomorphism between $\Galg_0^{F_0^2}\semi{F_0}$ and~$\wG^F$
constructed in the proof of Proposition~\ref{isogroupe}.

\end{remarque}

\section{Lusztig's conjecture for the Suzuki and the Ree groups}\label{sec3}

In this section, we will prove that the simple groups of
types~$B_2$,~$G_2$ and~$F_4$ with the Frobenius map defining the Suzuki
and Ree groups satisfy Lusztig's conjecture. For the explicit calculations in
this section, we will
use programs developed in~\textsc{Gap}~\cite{gap} by the author.
We will first recall some generalities.

\subsection{Generalized Springer correspondence}\label{secSP}
Let~$\Galg_0$ be a connected reductive group defined over the finite
field of~$q$ elements and~$F_0$ be the corresponding Frobenius map.
We set~$\mathcal{N}$ to be the~$\Galg_0$-classes of pairs~$(u,\phi)$,
where~$u$ is a unipotent element of~$\Galg_0$ and~$\phi$ an
irreducible character of the component group~$A_{\Galg_0}(u)$
defined in the proof of Proposition~\ref{classesext}. Two
pairs~$(u,\phi)$ and~$(v,\psi)$ are in the same~$\Galg_0$-class 
if there existe~$x\in\Galg_0$ such that~$v=x^{-1}ux$ and~$\psi(t)=\phi(xtx^{-1})$ 
for every~$t\in A_{\Galg_0}(v)$ (note that it is independent of the choice of~$x$).
In~\cite{Lustinter}, Lusztig associates to every pair~$(u,\phi)$ a
$4$-tuple~$(\Lalg,v,\psi,\rho)$, where~$\Lalg$ is a Levi subgroup
of~$\Galg_0$,~$v\in\Lalg$ is unipotent,~$\psi\in\Irr(A_{\Lalg}(v))$ and
$\rho\in\Irr(\Nn_{\Galg_0}(\Lalg)/\Lalg)$. If~$\Lalg=\Galg_0$ then the
pair is said {\it cuspidal}. Lusztig's construction defines a bijection
between~$\mathcal{N}$ and the set of $4$-tuples as above, such that
$(v,\psi)$ is cuspidal for~$\Lalg$. We write~$\rho_{u,\phi}$ for~$\rho$
if~$(u,\phi)$ corresponds to~$(\Lalg,v,\psi,\rho)$. The map~$(u,\phi)\mapsto\rho_{u,\phi}$ 
is the generalized Springer correspondence.

\subsection{Lusztig's conjecture}\label{secLC}

We use the same notation as in~\S\ref{secSP}. We
denote by~$\Talg_0$ a maximal rational torus of~$\Galg_0$ contained
in a rational Borel~$\Balg_0$ of~$\Galg_0$. We define the Weyl
group~$W=\Nn(\Talg_0)/\Talg_0$ of~$\Galg_0$. Since~$\Talg_0$ is
rational, it follows that~$F_0$ acts as an automorphism on~$W$;
we denote it with the same symbol. Let~$\rho\in\Irr(W)$ such
that~$\rho^{F_0}=\rho$; using Clifford theory we can extend~$\rho$
to~$\widetilde{W}=W\semi{F_0}$. We choose such
an extension~$\widetilde{\rho}$ and we define the \emph{uniform almost characters
of}~$\Galg_0^{F_0}$ as 
$$R_{\widetilde{\rho}}=\frac{1}{|W|}\sum_{w\in
W}\widetilde{\rho}(wF_0)R_w,
$$
where~$R_w$ is the generalized
Deligne-Lusztig character~$R_{\Talg_w}^{\Galg_0}(1_{T_w^{F_0}})$ associated
to the rational torus~$\Talg_w$ of~$\Galg_0$ corresponding to~$w$. For
definition and properties of these characters we refer for example
to~\cite[\S7]{Carter2}.

Let~$C$ be a rational unipotent class in~$\Galg_0$ and~$\mathcal{E}$
be an~$F_0$-stable irreducible local system on~$\Galg_0$ which
is~$\Galg_0$-equivariant for the conjugation action of~$\Galg_0$.
In~\cite{Lustsheave}, Lusztig associates to~$(C,\mathcal{E})$ an
irreducible $F_0$-stable perverse sheaf on~$\Galg_0$. We can attach to it
its characteristic function~$X_{C,\mathcal{E}}$ on~$\Galg_0^{F_0}$ (\cite[II.8.4]{Lustsheave}
for the definition), defined up to a scalar. 

Following~\cite[3.5]{shojiSC}, there is a correspondence between the
pairs~$(C,\mathcal{E})$ defined previously and the pairs~$(u,\phi)$
defined in~\S\ref{secSP} such that~$u$ is rational and~$\phi$
is~$F_0$-stable (the set of such pairs is denoted by~$\mathcal{N}^{F_0}$ in
the following); we denote by~$X_{u,\phi}$ the corresponding
characteristic function in this correspondence. We denote
by~$\mathcal{N}_0$ (resp.~$\mathcal{N}_0^{F_0}$) the pairs
of~$\mathcal{N}$ (resp.~$\mathcal{N}^{F_0}$) such that the 
associated $4$-tuple has the form~$(\Talg_0,1,1,\rho_{u,\phi})$. 
Note that the restriction of the map~$(u,\phi)\mapsto\rho_{u,\phi}$ 
to~$\mathcal{N}_0$ is the original Springer correspondence.

Lusztig conjectures that the uniform almost
characters~$R_{\widetilde{\rho}_{u,\phi}}$ and the characteristic
functions~$X_{u,\phi}$ with~$(u,\phi)\in\mathcal{N}_0^{F_0}$ coincide up
to a scalar. We propose to verify that this conjecture holds on
unipotent elements.

\begin{remarque} Lusztig's conjecture is more general. Lusztig defines
in~\cite[13.6]{Lust} the almost characters.
Lusztig's
conjecture asserts that these class functions coincide up to a scalar
with the characteristic functions of every $F_0$-stable character sheaf
on~$\Galg_0$. To compute the values of these functions, it is sufficient
to compute the values of the characteristic functions~$X_{C,\mathcal{E}}$ 
on the unipotent elements;
see~\cite[II.8.5]{Lustsheave}. 
\end{remarque}

\subsection{Lusztig's algorithm}\label{secAlgo} We keep the notation
of the  preceding sections. Our aim is to compute the values of the
characteristic functions~$X_{u,\phi}$ for~$(u,\phi)\in\mathcal{N}_0^{F_0}$
on the unipotent elements of~$\Galg_0^{F_0}$ and to compare them with
the uniform almost characters of~$\Galg_0^{F_0}$.

Let~$(u,\phi)\in\mathcal{N}_0^{F_0}$; we 
define~$\varphi_{u,\phi}:\Cl_{\Galg_0}(u)^{F_0}\rightarrow\overline{\Q}_
{\ell}$ as follows: using Theorem~\ref{gk}, if~$a$ is a representative
of an~$F_0$-class of~$A_{\Galg_0}(u)$, we denote by~$g_a$ the
corresponding class representative of~$\Galg_0^{F_0}$.
We set~$\varphi_{u,\phi}(g)=\widetilde{\phi}(aF_0)$ if~$g$ is
conjugate to~$g_a$ in~$\Galg_0^{F_0}$ and~$\varphi_{u,\phi}(g)=0$
if~$g\notin\Cl_{\Galg_0}(u)^{F_0}$; it depends on the choice of an extension~$\widetilde{\phi}$
of the~$F_0$-stable character~$\phi$
to~$A_{\Galg_0}(u)\semi{F_0}$ and hence it is only well-defined up to a non-zero scalar multiple. 
Note that the class function~$\varphi_{u,\phi}$
is up to a scalar the class function~$Y_{u,\phi}$ defined
in~\cite[V.24.2.3]{Lustsheave}.

We now recall Lusztig's algorithm. We first choose an order
on~$\mathcal{N}_0^{F_0}$ extending the natural order on unipotent
classes of~$\Galg_0^{F_0}$ defined by~$u\leq v$ when~$\Cl_{\Galg_0}(u)$
is contained in the Zariski closure of~$\Cl_{\Galg_0}(v)$. Following
Lusztig, we write~$(u,\phi)\sim(v,\psi)$ when $u$ is conjugate
to~$v$. If~$u$ is a unipotent element of~$\Galg_0$, we denote
by~$\mathcal{B}_u$ the variety of Borel subgroups of~$\Galg_0$
containing~$u$. We set~$d_u=\dim(\mathcal{B}_u)$. We remark that
if~$u\leq v$ then~$d_v\leq d_u$. In~\cite{Lustsheave}, Lusztig shows
that the funtions~$Y_{u,\phi}$ (with~$(u,\phi)\in\mathcal{N}_0^{F_0}$)
form a basis of the uniform functions on the set of unipotent elements
of~$\Galg_0^{F_0}$. Note that the function~$X_{u,\phi}$ defined
in~\cite[V.24.2.8]{Lustsheave} is the characteristic function
of the sheaf divided by~$q^{-(a_0+r)/2}$; in our case, we
have~$a_0+r=\dim(\Cen_{\Galg_0}(u))-\rank(\Galg_0)=2d_u$;
see~\cite[2.7]{Sp}. With our notation, we have
\begin{equation}\label{eq1}
q^{-d_u}X_{u,\phi}=\sum_{(v,\psi)\in\mathcal{N}_0^{F_0}}p_{(v,\psi),(u
,\phi)}Y_{v,\psi}.
\end{equation} 
We set~$P$ to be the matrix with
coefficients~$p_{(v,\psi),(u,\phi)}$. In~\cite[V.24]{Lustsheave},
Lusztig shows that the coefficients of~$P$ are rational, that the
diagonal elements of~$P$ are~$1$ and that~$P$ is an upper triangular
matrix. Moreover, if~$(u,\phi)\sim(v,\psi)$ but~$(u,\phi)\neq(v,\psi)$
then~$p_{(v,\psi),(u,\phi)}=0$. We give a way to compute
the matrix~$P$. We introduce the nonsingular bilinear form
$$(f,f')=\sum_{g\in X}f(g)\widetilde{f'}(g),$$ where~$X$ denotes
the set of unipotent elements of~$\Galg_0^{F_0}$
and~$\widetilde{f'}$ is defined in~\cite[V.24.2.12]{Lustsheave}. We
set $$w_{(u,\phi),(v,\psi)}=\frac{|\Galg_0^{F_0}|}{|W|}\sum_{w\in
W}\frac{\widetilde{\rho}_{u,\phi}(wF_0)\widetilde{\rho}_{v,\psi}(wF_0)}
{|\Talg_w^{F_0}|q^{d_u+d_v}},$$ where~$\widetilde{\rho}_{u,\phi}$
is an extension of the~$F_0$-stable character~$\rho_{u,\phi}$ given
in~\S\ref{secSP} to the group~$W\semi{F_0}$ and~$\Talg_w$ is the
rational torus of~$\Galg_0$ corresponding to~$w$.

We set~$\Omega$ to be the matrix with coefficients~$w_{(u,\phi),(v,\psi)}$.
Lusztig proves that the matrix~$P$ defined as above is the unique
solution of $$^tP\Lambda P=\Omega,$$ where~$\Lambda$ is a block diagonal
matrix, which is uniquely determined from the equation when we choose
the dimension of the blocks compatible with the relation~$\sim$.

\subsection{The Suzuki groups}

Let~$\Galg_0$ be a simple algebraic group of type~$B_2$ defined over
the algebraic closure of the field of~$2$ elements. Let~$n$ be a
non-negative integer; we set~$q=2^n\sqrt{2}$ and~$F_0=\s\circ F_{2^n}$,
where~$F_{2^n}$ is the standard Frobenius map on~$\Galg_0$ over~$\F_{2^n}$ and~$\s$
is the bijective endomorphism arising from the symmetry of the Dynkin diagram,
see~\cite[\S12.3]{Carter}. The Suzuki group of parameter~$q^2$ is defined
as~$\Galg_0^{F_0}$; it has order~$q^4(q^2-1)(q^4+1)$. The Weyl group of~$\Galg_0$ is the dihedral
group of order~$8$. We denote by~$w_a$ and~$w_b$ the two generators
of~$W$ corresponding to the roots~$a$ and~$b$ respectively. We denote
by~$x_r(t)$ (for a root $r$ and~$t\in\overline{\F}_2$) the Chevalley
generators of~$\Galg_0$ and we set~$x_a=x_a(1)$ and~$x_b=x_b(1)$. By
convention, we choose~$b$ for the long root. The group~$W$ has~$5$
irreducible characters. We denote by~$\varepsilon$ the sign character of~$W$
and by~$\chi$ the irreducible character of degree~$2$. The two other
linear characters are denoted by~$\varphi_1$ and~$\varphi_2$, the choice is such
that~$\varphi_1(w_a)=-1$.

In Table~\ref{SPB2}, we recall the generalized Springer
correspondence for a simple group of type~$B_2$
in characteristic~$2$;~see~\cite[IV.1]{Sp}. We
set~$u_1=1$,~$u_2=x_b$,~$u_3=x_a$, $u_4=x_{a+b}x_{2a+b}$
and~$u_5=x_ax_bx_{a+b}$ as representatives for the unipotent
classes of~$\Galg_0$. The~$F_0$-stable classes of~$\Galg_0$ have
representatives~$u_1$,~$u_4$ and~$u_5$. The class~$\Cl(u_5)^{F_0}$
splits into two classes of~$\Galg_0^{F_0}$; representatives of these
classes are~$x_ax_bx_{a+b}$ and~$x_ax_bx_{2a+b}$ and are denoted by~$\rho$
and~$\rho^{-1}$ respectively in~\cite{Suzuki}.

\begin{table}
\footnotesize
$$\renewcommand{\arraystretch}{1.5}\begin{array}{c|cccc}
u &A_{\Galg_0}(u)&\phi&\rho_{u,\phi}&d_u \\\hline
u_1&1&1&\varepsilon&4 \\
u_2&1&1&\varphi_1&2 \\
u_3&1&1&\varphi_2&2 \\
u_4&1&1&\chi&1 \\
u_5&\Z/2\Z&1&1_W&0 \\
\end{array}$$
\normalsize
\caption{Generalized Springer correspondence for type~$B_2$.}\label{SPB2}
\end{table}

\subsubsection{Uniform almost characters}

The Weyl group~$W$ of~$\Galg_0$ has three~$F_0$-stable
characters~$1_W$,~$\varepsilon$ and~$\chi$. We choose 
extensions~$\widetilde{1}_W$,~$\widetilde{\varepsilon}$
and~$\widetilde{\chi}$ to~$W\semi{F_0}$. The values of these characters
appear in Table~\ref{WB2}. \begin{table}\footnotesize
$$\renewcommand{\arraystretch}{1.5}\begin{array}{c|ccc}
 &F_0&w_aF_0&w_aw_bw_aF_0 \\\hline
\widetilde{\varepsilon}&1&-1&-1\\
\widetilde{\chi}&0&-\sqrt{2}&\sqrt{2}\\
\widetilde{1}_W&1&1&1\\
\end{array}$$\normalsize
\caption{Values of extensions on~$WF_0$.}\label{WB2}
\end{table} 
The decomposition into irreducible components of the uniform
almost characters (computed using the definition of uniform almost
characters in~\S\ref{secLC}) of the Suzuki
groups are given in~\cite[Table 1]{GM}. We
have~$R_{\widetilde{1}_W}=1_{\Galg_0^{F_0}}$,~$R_{\widetilde{\varepsilon
}}=\St_{\Galg_0^{F_0}}$
and~$R_{\widetilde{\chi}}=\frac{1}{\sqrt{2}}(W_1+W_2)$, where~$W_1$
and~$W_2$ are the two cuspidal unipotent characters of the Suzuki group
(see~\cite{Suzuki}); note that the results there are correct, however, one has 
to use our choices for~$\widetilde{\chi}$ given in Table~\ref{WB2} contrary to what they claim in~\cite[Table 1]{GM}. 
Thus, using the character table of the Suzuki group given in~\cite{Suzuki}, we deduce

\begin{proposition}\label{alb2}
The values of the uniform almost characters of the Suzuki groups 
on the unipotent elements of~$\Galg_0^{F_0}$ are given in
Table~\ref{almostB2}.
\end{proposition}

\begin{table}\footnotesize
$$\renewcommand{\arraystretch}{1.5}\begin{array}{c|cccc}
 &u_1&u_4&\rho&\rho^{-1}\\\hline
R_{\widetilde{\varepsilon}}&q^4&0&0&0\\
R_{\widetilde{\chi}}&q(q^2-1)&-q&0&0\\
R_{\widetilde{1}_W}&1&1&1&1\\
\end{array}$$\normalsize
\caption{Values of the uniform almost characters of the Suzuki group.}\label{almostB2}
\end{table}

\subsubsection{Lusztig's algorithm for the Suzuki groups}

The order on~$\mathcal{N}_0^{F_0}$ is~$(u_1,1)\leq(u_4,1)\leq (u_5,1)$.
For the construction of the functions~$\varphi_{u,\phi}$, 
we need to fix an extension of~$\phi$ to~$A_{\Galg_0}(u)\semi{F_0}$;
we choose the trivial character for~$A_{\Galg_0}(u_1)\semi{F_0}$
and for~$A_{\Galg_0}(u_3)\semi{F_0}$, and the non-trivial character
for~$A_{\Galg_0}(u_4)\semi{F_0}$. We then give the resulting
functions in Table~\ref{YB2}. Recall that~$|\Talg_1^{F_0}|=q^2-1$,
$|\Talg_{w_a}^{F_0}|=q^2-\sqrt{2}q+1$
and~$|\Talg_{w_aw_bww_a}^{F_0}|=q^2+\sqrt{2}q+1$.

\begin{table}\footnotesize
$$\renewcommand{\arraystretch}{1.5}\begin{array}{c|cccc}
 &u_1&u_4&\rho&\rho^{-1}\\\hline
\varphi_{u_1,1}&1&0&0&0\\
\varphi_{u_4,1}&0&-1&0&0\\
\varphi_{u_5,1}&0&0&1&1\\
\end{array}$$\normalsize
\caption{Values of the functions~$\varphi_{u,\phi}$ for the Suzuki group.}\label{YB2}
\end{table}

\begin{theoreme}When~$\Galg_0^{F_0}$ is the Suzuki group of 
parameter~$q^2$, then the matrix~$\Omega$ defined in~\S\ref{secAlgo} is
$$\Omega=\begin{bmatrix}
1&q^2-1&1\\
q^2-1&q^6-q^2&-q^6+q^4\\
1&-q^6+q^4&q^8
\end{bmatrix}
$$
The matrices~$P$ and~$\Lambda$ resulting from Lusztig's
algorithm described in~\S\ref{secAlgo} are
$$P=\begin{bmatrix}
1&q^2-1&1\\
0&1&-1\\
0&0&1
\end{bmatrix}\quad\textrm{and}\quad
\Lambda=\begin{bmatrix}
1&0&0\\
0&q^6-q^4+q^2-1&0\\
0&0&q^8-q^6+q^4-q^2
\end{bmatrix}$$ 
Moreover the characteristic functions~$X_{u,\phi}$ with~$(u,\phi)\in\mathcal{N}_0^{F_0}$ 
computed using Relation~(\ref{eq1}) in~\S\ref{secAlgo} coincide with the uniform almost
characters~$R_{\widetilde{\rho}_{u,\phi}}$ on the unipotent elements
of~$\Galg_0^{F_0}$.

\end{theoreme}

\subsection{The Ree groups of type~$G_2$}

Let~$\Galg_0$ be a simple algebraic group of type~$G_2$ defined
over~$\overline{\F}_3$; we set~$q=3^n\sqrt{3}$ for some non-negative
integer~$n$ and~$F_0=\s\circ F_{3^n}$, where~$\s$ is described
in~\cite[\S12.4]{Carter} and~$F_{3^n}$ is the standard Frobenius map over~$\F_{3^n}$.
The finite group~$\Galg_0^{F_0}$ is the
Ree group of type~$G_2$ with parameter~$q^2$ which has order~$q^6(q^2-1)(q^6+1)$.
The Weyl group
of~$\Galg_0$ is the dihedral group of order~$12$. We denote by~$w_a$
and~$w_b$ the reflections of~$W$ corresponding to the roots~$a$ and~$b$
respectively (here~$b$ is the long root). As before, we denote
by~$x_r(t)$ the Chevalley generators of~$\Galg_0$. The group~$W$ has~$6$
irreducible characters. We denote by~$1_W$ the trivial character
and by~$\varepsilon$ the sign of~$W$. There are two other linear
characters; we denote by~$\epsilon_a$ (resp.~$\epsilon_b$) the one which
satisfies~$\epsilon_a(w_a)=-1$ (resp.~$\epsilon_b(w_b)=-1$). We denote
by~$\theta'$ and~$\theta''$ the two characters of degree~$2$ such
that~$\theta'(w_aw_b)=1$.

In Table~\ref{SPG2}, we recall the generalized Springer correspondence
for a simple group of type~$G_2$ in characteristic~$3$;~see~\cite{spSP}.
To obtain the correspondence for~$\theta'$ and~$\theta''$ with
our notations, we use that~$b_{\theta'}=1$ and~$b_{\theta''}=2$,
where~$b_{\chi}$ denotes the smallest positive integer~$k$ such
that~$\chi$ occurs in the representation of~$W$ on the space of
homogeneous polynomials of degree~$k$ on the natural complex
reflection representation of~$W$. The representatives of unipotent classes
of~$\Galg_0$ are~$u_1=1$,~$u_2=x_{2a+b}(1)$,~$u_3=x_{3a+2b}(1)$,
$u_4=x_{2a+b}(1)x_{3a+2b}(1)$, $u_5=x_{a+b}(1)x_{3a+b}(1)$
and~$u_6=x_a(1)x_b(1)$. The~$F_0$-stable classes of~$\Galg_0$ have
representatives~$u_1$,~$u_4$,~$u_5$ and~$u_6$. The sets~$\Cl(u_5)^{F_0}$
and~$\Cl(u_6)^{F_0}$ split into two and three classes of~$\Galg_0^{F_0}$
respectively. The representatives are denoted by~$T$,~$T^{-1}$ and~$Y$,
$YT$, $YT^{-1}$ following the notation of~\cite{Ward}.

\begin{table}\footnotesize
$$\renewcommand{\arraystretch}{1.5}\begin{array}{c|cccc}
u &A_{\Galg_0}(u)&\phi&\rho_{u,\phi}&d_u \\\hline
u_1&1&1&\varepsilon&6 \\
u_2&1&1&\epsilon_a&3 \\
u_3&1&1&\epsilon_b&3 \\
u_4&1&1&\theta''&2 \\
u_5&\Z/2\Z&1&\theta'&1 \\
u_6&\Z/3\Z&1&1_W&0 \\
\end{array}$$
\normalsize
\caption{Generalized Springer correspondence for type~$G_2$.}\label{SPG2}
\end{table}

\subsubsection{Uniform almost characters}

The $F_0$-stable characters of~$W$ are~$1_W$,~$\varepsilon$,~$\theta'$ and~$\theta''$. 
In Table~\ref{WG2}, we give the values of an extension of one of 
their extensions on the coset~$W F_0$. To simplify the notation, 
we set~$w_1=w_a$, $w_2=w_aw_bw_a$ and~$w_3=w_aw_bw_aw_bw_a$.
\begin{table}\footnotesize
$$\renewcommand{\arraystretch}{1.5}\begin{array}{c|cccc}
 &F_0&w_1F_0&w_2F_0&w_3 F_0 \\\hline
\widetilde{\varepsilon}&1&-1&-1&-1\\
\widetilde{\theta}''&0&1&-2&1\\
\widetilde{\theta}'&0&-\sqrt{3}&0&\sqrt{3}\\
\widetilde{1}_W&1&1&1&1\\
\end{array}$$\normalsize
\caption{Values of extensions on~$WF_0$.}\label{WG2}
\end{table} 
In~\cite[Table 1]{GM}, the decomposition of the uniform almost 
characters of the Ree groups are given. We use the notation of~\cite{Ward} 
for the unipotent characters of~$\Galg_0^{F_0}$.  We have:
$$\begin{array}{lcl}
R_{\widetilde{1}_W}&=&1_{\Galg_0^{F_0}},\\
R_{\widetilde{\varepsilon}}&=&\St_{\Galg_0^{F_0}},\\
R_{\widetilde{\theta}'}&=&\frac{1}{2\sqrt{3}}(\xi_5+\xi_6+\xi_7+\xi_8+2\xi_9+2\xi_{10}),\\
R_{\widetilde{\theta}''}&=&\frac{1}{2}(\xi_5-\xi_6+\xi_7-\xi_8).\\
\end{array}$$ Using the character table of~$\Galg_0^{F_0}$ given
in~\cite{Ward}, we deduce

\begin{proposition}\label{alg2}The values of the uniform almost
characters on the unipotent elements of~$\Galg_0^{F_0}$ are given in Table~\ref{almostG2}.
\end{proposition}

\begin{table} 
\footnotesize
$$\renewcommand{\arraystretch}{1.5}\begin{array}{c|ccccccc}
 &u_1&u_4&T&T^{-1}&Y&YT&YT^{-1}\\\hline
R_{\widetilde{\varepsilon}}&q^6&0&0&0&0&0&0\\
R_{\widetilde{\theta}''}&q^2(q^2-1)&-q^2&0&0&0&0&0\\
R_{\widetilde{\theta}'}&q(q^4-1)&-q&-q&-q&0&0&0\\
R_{\widetilde{1}_W}&1&1&1&1&1&1&1\\
\end{array}$$
\normalsize
\caption{Values of the uniform almost characters of the Ree groups of type~$G_2$.}\label{almostG2}
\end{table}

\subsubsection{Lusztig's algorithm for the Ree groups of type~$G_2$}

We first recall the orders of the finite tori of~$\Galg_0^{F_0}$
corresponding to~$w\in W$; we have~$|\Talg_1^{F_0}|=q^2-1$,
$|\Talg_{w_1}^{F_0}|=q^2-\sqrt{3}q+1$, $|\Talg_{w_2}^{F_0}|=q^2+1$
and $|\Talg_{w_3}^{F_0}|=q^2+\sqrt{3}q+1$. We choose the order
$(u_1,1)\leq(u_4,1)\leq(u_5,1)\leq(u_6,1)$ on~$\mathcal{N}_0^{F_0}$.
In order to define the~$\varphi_{u,\phi}$ we choose as extensions
the trivial character of~$A_{\Galg_0}(u_1)\semi{F_0}$,
of~$A_{\Galg_0}(u_4)\semi{F_0}$ and of~$A_{\Galg_0}(u_6)\semi{F_0}$.
For~$A_{\Galg_0}(u_5)\semi{F_0}$, we choose the extension of the trivial
character of~$A_{\Galg_0}(u_5)$ which is not the trivial character
of~$A_{\Galg_0}(u_5)\semi{F_0}$. In Table~\ref{YG2}, we give the values
of the corresponding functions.

\begin{table}\footnotesize
$$\renewcommand{\arraystretch}{1.5}\begin{array}{c|ccccccc}
 &u_1&u_4&T&T^{-1}&Y&YT&YT^{-1}\\\hline
\varphi_{u_1,1}&1&0&0&0&0&0&0\\
\varphi_{u_4,1}&0&1&0&0&0&0&0\\
\varphi_{u_5,1}&0&0&-1&-1&0&0&0\\
\varphi_{u_6,1}&0&0&0&0&1&1&1\\
\end{array}$$\normalsize
\caption{Values of the functions~$\varphi_{u,\phi}$ for the Ree groups of type~$G_2$.}\label{YG2}
\end{table}

\begin{theoreme}
When~$\Galg_0^{F_0}$ is the Ree group of type~$G_2$ with
parameter~$q^2$, then the matrix~$\Omega$ defined in~\S\ref{secAlgo} is
$$\Omega=\begin{bmatrix}
1&1-q^2&q^4-1&1\\
1-q^2&q^8-q^6-q^4-q^2&-q^8+q^4&q^8-q^6\\
q^4-1&-q^8+q^4&q^{10}+q^8-q^6-q^4&-q^{10}+q^6\\
1&q^8-q^6&-q^{10}+q^{6}&q^{12}
\end{bmatrix}
$$
The matrices~$P$ and~$\Lambda$ resulting from Lusztig's algorithm
described in~\S\ref{secAlgo} are 
$$P=\begin{bmatrix}
1&1-q^2&q^4-1&1\\
0&1&-1&1\\
0&0&1&-1\\
0&0&0&1
\end{bmatrix}\quad\textrm{and}\quad
\Lambda=\begin{bmatrix}
1&0&0&0\\
0&p_1(q)&0&0\\
0&0&q^2p_1(q)&0\\
0&0&0&q^4p_1(q)
\end{bmatrix},$$
where~$p_1(q)=(q^4-1)(q^4-q^2+1)$. 
Moreover for every~$(u,\phi)\in\mathcal{N}_0^{F_0}$ the
resulting characteristic function~$X_{u,\phi}$
coincides with the uniform almost
character~$R_{\widetilde{\rho}_{u,\phi}}$ on the unipotent elements
of~$\Galg_0^{F_0}$.
\end{theoreme}

\subsection{The Ree groups of type~$F_4$}

Let~$\Galg_0$ be a simple algebraic group of type~$F_4$ defined
over~$\overline{\F}_2$; we set~$q=2^n\sqrt{2}$ for some non-negative
integer~$n$ and~$F_0=\s\circ F_{2^n}$, where~$\s$ is described
in~\cite[\S12.3]{Carter}. The finite group~$\Galg_0^{F_0}$ is the
Ree group of type~$F_4$ with parameter~$q^2$ which has 
order~$q^{24}(q^2-1)(q^6+1)(q^8-1)(q^{12}+1)$.
The Weyl group of~$\Galg_0$ has~$1152$ elements and has~$25$ irreducible
characters. The character table of~$W$ was computed by Kondo
in~\cite{Kondo}; we use his notation. 

The generalized Springer correspondence of a simple algebraic
group of type~$F_4$ in characteristic~$2$ is determined
in~\cite[p.330]{spSP}. Remark that the characters denoted
by~$\chi_4$,~$\chi_{4,1}$, $\chi_{4,2}$, $\chi_{4,3}$, 
$\chi_{4,4}$,~$\chi_{12}$ and~$\chi_{16}$ in~\cite{spSP}
are those denoted by~$\chi_{4,1}$,~$\chi_{4,2}$, $\chi_{4,3}$, 
$\chi_{4,4}$, $\chi_{4,5}$,~$\chi_{12,1}$
and~$\chi_{16,1}$ in~\cite{Kondo} respectively. The~$20$ unipotent classes
of~$\Galg_0$ were classified by Shinoda and we use the notation
of~\cite[Table I]{Shinoda}. The representatives~$x_i$ are defined
as products of Chevalley generators in~\cite[p.139]{Shinoda2}. We
recall in Table~\ref{SPF4} the Springer correspondence of~$\Galg_0$. We
denote by~$\epsilon$ the non-trivial character of~$\Z/2\Z$, by~$\theta$
the irreducible character of~$\mathfrak{S}_3$ of degree~$2$, and by~$\epsilon'$
and~$\epsilon''$ the two linear characters of~$D_8$ distinct from the
trivial and the sign characters. There are~$10$ unipotent classes
of~$\Galg_0$ which are~$F_0$-stable. They have as representatives the
elements~$x_0$, $x_3$, $x_4$, $x_9$, $x_{15}$, $x_{16}$, $x_{17}$,
$x_{24}$, $x_{29}$ and~$x_{31}$. In~\cite{Shinoda}, Shinoda shows
that~$F_0$ acts trivially on~$A_{\Galg_0}(x_9)$, $A_{\Galg_0}(x_{29})$,
$A_{\Galg_0}(x_{17})$ and~$A_{\Galg_0}(x_{31})$. 
The action of~$F_0$ on~$A_{\Galg_0}(x_{24})$ is not trivial; 
this group has three~$F_0$-classes. Note that~$\epsilon'^{F_0}=\epsilon''$.
This permits to parametrize the conjugacy classes of~$\Galg_0^{F_0}$. We
denote by~$u_i$ ($0\leq i\leq 18$) a system of representatives as in~\cite[Table 2]{Shinoda}.

\begin{table}
\footnotesize
$$\begin{array}{cc}\renewcommand{\arraystretch}{1.5}\begin{array}{c|cccc}
u &A_{\Galg_0}(u)&\phi&\rho_{u,\phi}&d_u \\\hline
x_0&1&1&\chi_{1,4}&24 \\
x_1&1&1&\chi_{2,4}&16 \\
x_2&1&1&\chi_{2,2}&16 \\
x_3&1&1&\chi_{4,5}&13 \\
x_4&1&1&\chi_{9,4}&10 \\
x_5&\Z/2\Z&1&\chi_{8,4}&9 \\
&&\epsilon&\chi_{1,2}&\\
x_7&\Z/2\Z&1&\chi_{8,2}&9\\
&&\epsilon&\chi_{1,3}&\\
x_9&\Z/2\Z&1&\chi_{4,1}&8\\
x_{11}&1&1&\chi_{4,3}&7\\
x_{12}&1&1&\chi_{4,4}&7 \\
x_{13}&1&1&\chi_{9,2}&6 \\
\end{array}&\qquad
\renewcommand{\arraystretch}{1.5}\begin{array}{c|cccc}
u &A_{\Galg_0}(u)&\phi&\rho_{u,\phi}&d_u \\\hline
x_{14}&1&1&\chi_{9,3}&6 \\
x_{15}&1&1&\chi_{6,1}&6 \\
x_{16}&1&1&\chi_{16,1}&5 \\
x_{17}&\mathfrak{S}_3&1&\chi_{12,1}&4 \\
&&\theta&\chi_{6,2}\\
x_{20}&\Z/2\Z&1&\chi_{8,3}&3 \\
x_{22}&\Z/2\Z&1&\chi_{8,1}&3 \\
x_{24}&D_8&1&\chi_{9,1}&2\\
&&\epsilon'&\chi_{2,1}\\
&&\epsilon''&\chi_{2,3}\\
x_{29}&\Z/2\Z&1&\chi_{4,2}&1 \\
x_{31}&\Z/4\Z&1&\chi_{1,1}&0 \\
\end{array}
\end{array}
$$\normalsize
\caption{Springer correspondence for type~$F_4$.}\label{SPF4}
\end{table}

\subsubsection{Uniform almost characters}

The group~$W$ has eleven $F_0$-stable characters:
$\chi_{1,1}$,~$\chi_{1,4}$, $\chi_{4,1}$, $\chi_{4,2}$, $\chi_{4,5}$,
$\chi_{6,1}$, $\chi_{6,2}$, $\chi_{9,1}$, $\chi_{9,4}$, $\chi_{12,1}$
and~$\chi_{16,1}$. For each of these characters, we choose an extension 
to~$W\semi{F_0}$. We give the values of these extensions on the
coset~$WF_0$ in~Table~\ref{WF4}. We take the elements~$w_i$
($1\leq i \leq 11 $) defined in~\cite[p.8]{Shinoda} for a system of
representatives of the~$F_0$-classes of~$W$.

\begin{table}
\scriptsize
$$\renewcommand{\arraystretch}{1.3}\begin{array}{c|ccccccccccc}
&w_1F_0&w_2F_0&w_3F_0&w_4F_0&w_5F_0&w_6F_0&w_7F_0&w_8F_0&w_9F_0&w_{10}F_0&w_{11}F_0\\\hline
\widetilde{\chi}_{1,4}&1&-1& -1& -1& 1& 1& 1& 1& 1& 1& 1\\
\widetilde{\chi}_{4,5}&0&0&-\sqrt{2}&\sqrt{2}&0&2\sqrt{2}&-2\sqrt{2}& 0& 0&\sqrt{2}& -\sqrt{2}\\
\widetilde{\chi}_{9,4}& 1&1& -1& -1& -1& 3& 3&  -3& 0& 0& 0\\
\widetilde{\chi}_{4,1}& 2&0& 0& 0& 2& 2& 2& 2& -1& -1& -1\\
\widetilde{\chi}_{6,1}& 0&0& 0& 0& 2& -2& -2& -4& -1& 1& 1\\
\widetilde{\chi}_{16,1}&0&0& 0& 0& 0& 4\sqrt{2}& -4\sqrt{2}& 0& 0& -\sqrt{2}&\sqrt{2}\\
\widetilde{\chi}_{6,2}&-2&0& 0& 0& 0& 4& 4& 2& -1& 1& 1\\
\widetilde{\chi}_{12,1}& 2&0& 0& 0& -2& -2& -2& 2& -1& 1& 1\\
\widetilde{\chi}_{9,1}& 1&-1& 1& 1& -1& 3& 3& -3& 0& 0& 0\\
\widetilde{\chi}_{4,2}&0&0& -\sqrt{2}& \sqrt{2}& 0& -2\sqrt{2}& 2\sqrt{2}&0& 0& -\sqrt{2}& \sqrt{2}\\
\widetilde{\chi}_{4,1}& 1&1& 1& 1& 1& 1& 1& 1& 1& 1& 1
\end{array}$$
\normalsize
\caption{Values of extensions on~$WF_0$.}\label{WF4}
\end{table}

The unipotent characters of~$\Galg_0^{F_0}$ were computed by Malle
in~\cite{Malle2F4}; in the following, we will use his notation. The
decomposition of the uniform almost characters of~$\Galg_0^{F_0}$ in
unipotent components is given in~\cite[Table 1]{GM}. When~$\chi_{i,j}$
is an~$F_0$-stable irreducible character of~$W$, we write~$R_{i,j}$ 
for~$R_{\widetilde{\chi}_{i,j}}$ to simplify. We recall:
$$\begin{array}{lcl}
R_{1,4}&=&\chi_4,\\
R_{4,5}&=&\frac{1}{\sqrt{2}}(\chi_7+\chi_8),\\
R_{9,4}&=&\chi_3,\\
R_{4,1}&=&\frac{1}{4}(\chi_9+\chi_{10}+2\chi_{11}-\chi_{12}-\chi_{13}-2\chi_{14}-\chi_{15}-\chi_{16}-\chi_{17}-\chi_{18}),\\
R_{6,1}&=&\frac{1}{3}(\chi_{12}+\chi_{13}-\chi_{14}+\chi_{19}+\chi_{20}-2\chi_{21}),\\
R_{16,1}&=&\frac{1}{2\sqrt{2}}(\chi_9-\chi_{10}-\chi_{12}+\chi_{13}+\chi_{15}+\chi_{16}-\chi_{17}-\chi_{18}),\\
R_{6,2}&=&\frac{1}{6}(-3\chi_9-3\chi_{10}-\chi_{12}-\chi_{13}-2\chi_{14}+2\chi_{19}+2\chi_{20}+2\chi_{21}),\\
R_{12,1}&=&\frac{1}{12}(3\chi_9+3\chi_{10}+6\chi_{11}+\chi_{12}+\chi_{13}
+2\chi_{14}+3\chi_{15}+3\chi_{16}+3\chi_{17}\\&&+3\chi_{18}+4\chi_{19}+4
\chi_{20}+4\chi_{21}),\\
R_{9,1}&=&\chi_2,\\
R_{4,2}&=&\frac{1}{\sqrt{2}}(\chi_5+\chi_6),\\
R_{4,1}&=&\chi_1.
\end{array}$$
Using the table in~\cite{Malle2F4}, we deduce

\begin{proposition}The
values on the unipotent elements of the uniform almost 
characters of Ree groups of type~$F_4$ are given in Table~\ref{almostF4}.
\end{proposition}

\begin{sidewaystable}
\scriptsize
$$
\begin{array}{c}
\renewcommand{\arraystretch}{1.5}
\begin{array}{l|llllll}
&u_0&u_1&u_2&u_3&u_4&u_5\\\hline
R_{1,4}& q^{24}&0& 0& 0& 0& 0\\
R_{4,5}&q^{23}-q^{19}+q^{17}-q^{13}& -q^{13}& 0& 0& 0& 0\\
R_{9,4}& q^{22}-q^{20}+q^{16}-q^{12}+q^{10}& -q^{12}+q^{10}& q^{10}& 0& 0& 0\\
R_{4,1}&q^{16}+q^8& q^8& q^8& q^8& q^8& 0\\
R_{6,1}&q^{18}-q^{16}-q^8+q^6& -q^8+q^6& -q^8+q^6& 0& 0& q^6\\
R_{16,1}&-q^{19}+2q^{15}-q^{13}-q^{11}+2q^9-q^5& -q^{11}+2q^9-q^5& q^9-q^5& q^7-q^5& q^7-q^5& -q^5\\
R_{6,2}& -q^{18}+q^{16}-2q^{12}+q^8-q^6& -q^{12}+q^8-q^6& -q^6& -q^4& -q^4& 0\\
R_{12,1}&q^{20}+q^4& q^4& q^4& q^4& q^4& q^4\\
R_{9,1}&q^{14}-q^{12}+q^8-q^4+q^2& q^8-q^4+q^2& -q^4+q^2& q^6-q^4+q^2& q^6-q^4+q^2& -q^4+q^2\\
R_{4,2}& q^{11}-q^7+q^5-q& -q^7+q^5-q& q^5-q& q^3-q& q^3-q& -q\\
R_{4,1}& 1& 1& 1& 1& 1& 1\end{array}\\
\\
\\
\renewcommand{\arraystretch}{1.5}
\begin{array}{l|lllllllllllll}
&u_6&u_7&u_8&u_9&u_{10}&u_{11}&u_{12}&u_{13}&u_{14}&u_{15}&u_{16}&u_{17}&u_{18}\\ \hline
R_{1,4}& 0& 0& 0& 0& 0& 0& 0& 0& 0& 0&0&0&0\\
R_{4,5}& 0& 0& 0& 0& 0& 0& 0& 0& 0& 0&0&0&0\\
R_{9,4}& 0& 0& 0& 0& 0& 0& 0& 0& 0& 0&0&0&0\\
R_{4,1}& 0& 0& 0& 0& 0& 0& 0& 0& 0& 0&0&0&0\\
R_{6,1}& 0& 0& 0& 0& 0& 0& 0& 0& 0& 0&0&0&0\\
R_{16,1}&-q^5& 0& 0& 0& 0& 0& 0& 0& 0& 0&0&0&0\\
R_{6,2}&-q^4& -2q^4& 0& q^4& 0& 0& 0& 0& 0& 0& 0& 0& 0\\
R_{12,1}&q^4& q^4& q^4&q^4& 0& 0& 0& 0& 0& 0& 0& 0& 0\\\
R_{9,1}&-q^4+q^2& q^2& q^2& q^2& q^2& q^2& q^2& 0& 0& 0& 0& 0& 0\\
R_{4,2}&q^3-q& 2q^3-q& -q& -q^3-q& -q& -q& -q&-q&-q& 0& 0& 0& 0\\
R_{4,1}&1& 1& 1& 1& 1& 1& 1& 1& 1& 1&1&1&1
\end{array}
\end{array}$$
\normalsize
\caption{Values on unipotent elements of the uniform almost characters of Ree groups of type~$F_4$.\label{almostF4}
} 
\end{sidewaystable}

\subsubsection{Lusztig's algorithm for the Ree groups of type~$F_4$}
The order on~$\mathcal{N}_0^{F_0}$ is such that
\small
$$(1,4)\leq(4,5)\leq(9,4)\leq(4,1)\leq(6,1)\leq(16,1)\leq(6,2)\leq(12,1)\leq(9,1)\leq(4,2)\leq(1,1).$$
\normalsize
The orders of the finite tori corresponding to elements~$w\in W$ are
\footnotesize
$$\begin{array}{cc}
\renewcommand{\arraystretch}{1.4}
\begin{array}{l|l}
w&|\Talg_w^{F_0}|\\\hline
w_1&(q^2-1)^2\\
w_2&q^4-1\\
w_3&(q^2-1)(q^2-\sqrt{2}q+1)\\
w_4&(q^2-1)(q^2+\sqrt{2}q+1)\\
w_5&q^4+1\\
w_6&(q^2-\sqrt{2}q+1)^2\\
\end{array}&\qquad
\renewcommand{\arraystretch}{1.4}
\begin{array}{l|l}
w&|\Talg_w^{F_0}|\\\hline
w_7&(q^2+\sqrt{2}q+1)^2\\
w_8&(q^2+1)^2\\
w_9&q^4-q^2+1\\
w_{10}&q^4-\sqrt{2}q^3+q^2-\sqrt{2}q+1\\
w_{11}&q^4+\sqrt{2}q^3+q^2+\sqrt{2}q+1
\end{array}
\end{array}
$$\\
\normalsize

\noindent Using Table~\ref{WF4}, we can compute the matrix~$\Omega$
defined in~\S\ref{secAlgo}. We now give the selected
extension of~$\phi$ to~$A_{\Galg_0^{F_0}}(u)\semi{F_0}$
for~$(u,\phi)\in\mathcal{N}_0^{F_0}$. We choose the trivial character
of~$A_{\Galg_0^{F_0}}(u)\semi{F_0}$ for~$(x_0,1)$, $(x_4,1)$, $(x_9,1)$,
$(x_{15},1)$, $(x_{17},1)$, $(x_{24},1)$ and~$(x_{31},1)$. We choose the extension of the
trivial character which is not the trivial character for~$(x_3,1)$,
$(x_{16},1)$ and~$(x_{29},1)$. The group~$A_{\Galg_0}(x_{17})$
is isomorphic to~$\mathfrak{S}_3$ and~$F_0$ acts trivially on
it. Thus~$A_{\Galg_0}(x_{17})\semi{F_0}$ is a direct product; we
choose for~$\widetilde{\theta}$ the extension of~$\theta$ such
that~$\widetilde{\theta}(F_0)=-2$. We then can compute the corresponding functions~$\varphi_{u,\phi}$.

\begin{theoreme}
When~$\Galg_0^{F_0}$ is the Ree group of type~$F_4$ with parameter~$q^2$,
then the matrix~$P$ resulting from Lusztig's algorithm is
\footnotesize
$$P=\left[\renewcommand{\arraystretch}{1.5}\begin{array}{ccccccccccc}
1& p_1(q)& p_2(q)& p_3(q)& p_4(q)& p_5(q)&p_6(q)& p_7(q)& p_8(q)&p_1(q)&1\\
0& 1& p_9(q)& -1& p_9(q)& p_{10}(q)& p_{11}(q)&-1& p_{12}(q)& p_{13}(q)&-1 \\
0& 0& 1& 1& -p_9(q)& p_{14}(q)& -q^2& 1& -p_9(q)& p_{14}(q)& 1\\ 
0& 0& 0& 1& 0& p_{9}(q)& -1& 1& p_{15}(q)& p_9(q)& 1 \\
0& 0& 0& 0& 1& -1& 0& 1& -p_9(q)&-1& 1 \\
0& 0& 0& 0& 0& 1& 1& -1& p_9(q)& -p_9(q)& -1\\ 
0& 0& 0& 0& 0& 0& 1& 0& 0& -q^2&0\\
0& 0& 0& 0& 0& 0& 0& 1& 1& -1& 1 \\
0& 0& 0& 0& 0& 0& 0& 0& 1& -1& 1\\
0& 0& 0& 0& 0& 0& 0& 0& 0& 1& -1\\
0& 0& 0& 0& 0& 0& 0& 0& 0& 0& 1
\end{array}\right],$$
\normalsize
where
\footnotesize
$$\begin{array}{lcl}
p_1(q)&=&(q^2-1)(q^2+1)^2(q^4-q^+1)\\
p_2(q)&=&(q^8-q^4+1)(q^4-q^2+1)\\
p_3(q)&=&q^8+1\\
p_4(q)&=&(q^2-1)^2(q^4-q^3+q^2-q+1)(q^4+q^3+q^2+q+1)\\
p_5(q)&=&-(q^2-1)^2(q^2+1)^3(q^4-q^2+1)\\
p_6(q)&=&-q^2(q^{12}-q^{10}+2q^6-q^2+1)\\
p_7(q)&=& q^{16}+1\\
p_8(q)&=&(q^4-q^2+1)(q^8-q^4+1)\\
p_9(q)&=&q^2-1\\
p_{10}(q)&=&(q^2-1)(q^4-q^2-1)\\
p_{11}(q)&=&q^2(q^6-q^2+1)\\
p_{12(q)}&=& -q^6+q^2-1\\
p_{13}(q)&=&q^6-q^4+1\\
p_{14}(q)&=&q^4-1\\
p_{15}(q)&=&q^4-q^2+1
\end{array}$$
\normalsize
The matrix~$\Lambda$ is the diagonal matrix with entries 
\small
$$
\left[1,f_1(q),q^4f_2(q),q^{10}f_1(q),q^8f_3(q),q^{10}f_3(q),q^{12}f_3(q),q^{12}f_3(q),q^{16}f_3(q),q^{18}f_3(q),q^{20}f_3(q)\right],
$$
\normalsize
where
\small
$$\begin{array}{lcl}
f_1(q)&=&(q^2-1)(q^2+1)^2(q^4+1)(q^4-q^2+1)(q^8-q^4+1)\\
f_2(q)&=&(q^4-1)(q^4+1)^2(q^4-q^2+1)(q^8-q^4+1)\\
f_3(q)&=&(q^2-1)^2(q^2+1)^2(q^4+1)^2(q^4-q^2+1)(q^8-q^4+1)\\
\end{array}$$
\normalsize

Moreover, for every~$(u,\phi)\in\mathcal{N}_0^{F_0}$, the characteristic
function~$X_{u,\phi}$ obtained using this algorithm coincides with the
uniform almost character~$R_{\widetilde{\rho}_{u,\phi}}$ on the
unipotent elements of~$\Galg_0^{F_0}$.
\end{theoreme}

\section{The disconnected cases}\label{sec4}
We keep the same notation as in the preceding sections.

\subsection{Uniform almost characters}

Digne-Michel~\cite{DMnonconnexe} and Malle~\cite{Malle1} have
generalized the Deligne-Lusztig theory to disconnected groups.
We briefly recall the construction: let~$\widetilde{\Talg}'$
be a rational "torus" of~$\widetilde{\Galg}$ contained in
a~"Borel"~$\widetilde{\Balg}'$. We then have 
(see~\cite[Proposition 1.5]{DMnonconnexe})~$\widetilde{\Balg}'=\Ualg'\rtimes\widetilde{\Talg}'$,
where~$\Ualg'$ is the unipotent radical of~$\widetilde{\Balg}'^{\circ}$.
We introduce the variety~$\mathbf{Y}_{\Ualg'}=\{x\in\Galg\ |\
x^{-1}F(x)\in\Ualg'\}.$ The groups~$\Galg^F$ and~$\widetilde{\Talg}'^F$
act on~$\mathbf{Y}_{\Ualg'}$ and these actions commute. This induces
linear actions on~$H_c^i(\mathbf{Y}_{\Ualg'},\overline{\Q}_{\ell})$,
where~$\ell$ is a prime number distinct from the characteristic.
Let~$\theta\in\Irr(\widetilde{\Talg}'^F)$, we define the associated
generalized Deligne-Lusztig characters by (see~\cite[\S2]{DMnonconnexe})
$$R_{\widetilde{\Talg}'}^{\widetilde{\Galg}}(\theta)(g)=
\frac{1}{|\widetilde{\Talg}'^F|}\sum_{t\in\widetilde{\Talg}'^F}\theta(t^{-1})\Tr\left((g,t)|H_c^*(\mathbf{Y}_{\Ualg'},
\overline{\Q}_{\ell})\right)\quad\textrm{for }g\in\widetilde{\Galg}.$$
We now give a parametrization (up to conjugacy) of the rational "tori" of~$\widetilde{\Galg}$ containing~$\tau$.
We recall that the "tori" of~$\widetilde{\Galg}$ are conjugate in~$\Galg$; see~\cite{DMnonconnexe}. 
Let~$\widetilde{\Talg}'$ be a "torus" of~$\widetilde{\Galg}$ 
containing~$\tau$, then~$\widetilde{\Talg}'^{\circ}=\widetilde{\Talg}'\cap\Galg^{\circ}$ is a 
maximal torus of~$\Galg$
and we have~$\widetilde{\Talg}'=\widetilde{\Talg}'^{\circ}\semi{\tau}$.
We have

\begin{proposition}\label{tores}We fix a maximal rational torus~$\Talg_0$ of~$\Galg_0$ and 
we set~$\Talg=\Talg_0\times\Talg_0$ and~$\widetilde{\Talg}=\Talg\semi{\tau}$.
The~$\widetilde{\Galg}^F$-classes of rational "tori" of~$\widetilde{\Galg}$ are in correspondence with
the $\Galg_0^{F_0}$-classes of maximal rational tori of~$\Galg_0$. We fix a system~$W_{F_0}$ of representatives 
of~$F_0$-classes of~$W_0$. Let~$w\in W_{F_0}$ and~$n_w$ 
be a corresponding element in~$\Nn_{\Galg_0}(\Talg_0)$. 
We set $$\Talg_{0,w}=x\Talg_0 x^{-1}\quad\textrm{and}\quad
\widetilde{\Talg}_w=(x,x)\widetilde{\Talg}(x,x)^{-1},$$ where~$x$ is such that~$n_w=x^{-1}F_0(x)$.
Then the sets~$\{\Talg_{0,w}\ |\ w\in W_{F_0}\}$ and~$\{\widetilde{\Talg}_w\ |\ w\in W_{F_0}\}$ are respectively representative systems of~$\Galg_0^{F_0}$-classes of maximal rational tori of~$\Galg_0$ 
and~$\widetilde{\Galg}^F$-classes of rational "tori" of~$\widetilde{\Galg}$.
\end{proposition}

\begin{preuve}{}
Write~$X$ for the set of "tori" of~$\widetilde{\Galg}$. 
The connected group~$\Galg$ acts transitively on~$X$. 
We denote by~$F':X\rightarrow X$, $\widetilde{\Talg}'\mapsto F(\widetilde{\Talg}')$. 
We have~$F'(x\widetilde{\Talg}')=F(x)F'(\widetilde{\Talg}')$
and~$\Nn_{\Galg^{\circ}}(\Talg')$ is a closed subset. We can apply Theorem~\ref{gk}.
The "torus"~$\widetilde{\Talg}=\left(\Talg_0\times\Talg_0\right)\semi\tau$ is a rational "torus" of~$X$. We have
$$\Nn_{\Galg}(\widetilde{\Talg})=\Nn_{\Galg}(\Talg_0\times\Talg_0)\cap\Nn_{\Galg}(\tau)=\mu\left(\Nn_{\Galg_0}(\Talg_0)\right),$$
it follows that~$\mu$ is an isomorphism of algebraic groups between~$\Nn_{\Galg_0}(\Talg_0)$ 
and~$\Nn_{\Galg}(\widetilde{\Talg})$. Since
$\Cen_{\Galg}(\widetilde{\Talg})=\mu\left(\Talg_0\right)$, we deduce that
$$\Nn_{\Galg}(\widetilde{\Talg})^{\circ}=\Cen_{\Galg}(\widetilde{\Talg}).$$
Moreover, we have~$$\Nn_{\Galg}(\widetilde{\Talg})/\Cen_{\Galg}(\widetilde{\Talg})\simeq
\mu\left(\Nn_{\Galg_0}(\Talg_0)/\Talg_0\right),$$
where~$\mu$ is the map~$w\Talg_0\mapsto
\mu(w\Talg_0)$. 
We then have
%
%
$$\forall w\in\Nn_{\Galg_0}(\Talg_0)/\Talg_0,\quad  F\left(\mu(w)\right)=\mu\left(F_0(w)\right).$$
Let~$w\in \Nn_{\Galg_0}(\Talg_0)/\Talg_0$ and~$n_w$ be a representative of~$w$ in~$\Nn_{\Galg_0}(\Talg_0)$. 
Let~$x\in\Galg_0$ be such that~$n_w=x^{-1}F_0(x)$;
We have~$$(x,x)^{-1}F(x,x)=\left(x^{-1}F_0(x),x^{-1}F_0(x)\right)=(n_w,n_w)=\mu(n_w),$$
The result is then a consequence of~Theorem~\ref{gk}.

\end{preuve}

\begin{remarque} 
We have to find an analogue of Proposition~1.40
of~\cite{DMnonconnexe} in our special case (with explicit representatives for the classes
of rational~"tori"). As a consequence of Proposition~\ref{tores},
we obtain Proposition~1.38 of~\cite{DMnonconnexe}, which says that
every~$\widetilde{\Galg}^F$-class of~"torus"~of~$\widetilde{\Galg}$ has
a representative that contains~$\tau$.
\end{remarque}
The Weyl group of~$\Galg$ is~$W=\Nn_{\Galg}(\Talg)/\Talg$, where~$\Talg=\Talg_0\times\Talg_0$. 
Following Malle in~\cite{Malle1}, we give a definition of the uniform almost characters of~$\widetilde{\Galg}^F$:
let~$\rho$ be an irreducible character of~$W^{\tau}$ and~$\widetilde{\rho}$ 
be an extension of~$\rho$ to~$\widetilde{W}={W^{\tau}}\semi{F}$. 
We define the uniform almost character corresponding to~$\widetilde{\rho}$ as the restriction of the class function onto the coset~$\Galg^F\tau$:
(we adapt here the definition~\cite[Definition 2]{Malle1} given in the case where~$F$ acts trivially on the Weyl group):
$$R_{\widetilde{\rho}}=\frac{1}{|W^{\tau}|}\sum_{w\in W^{\tau}}
\widetilde{\rho}(wF)R_{\widetilde{\Talg}_w}^{\widetilde{\Galg}}(1_{\widetilde{\Talg}_w^F}).$$
Note that~$W^{\tau}$ is isomorphic to~$W_0$ and the operation of~$F$ on~$W^{\tau}$ is the same as the one of~$F_0$ on~$W_0$.\\

In order to compute the decomposition of 
characters~$R_{\widetilde{\Talg}_w}^{\widetilde{\Galg}}(1_{\widetilde{
\Talg}_w^F})$ in irreducible constituents, we recall Proposition~4.8 of~\cite{DMnonconnexe}:

\begin{proposition}\label{propDM} Let~$\widetilde{\Talg}'$ and~$\widetilde{\Talg}''$
be two rational~"tori" containing~$\tau$; suppose 
that~$\theta'\in\Irr\left(\widetilde{\Talg}'^F\right)$
and~$\theta''\in\Irr\left(\widetilde{\Talg}''^F\right)$ satisfy
$\theta_0'=\Res_{\widetilde{\Talg}'^{\circ F}}^{\widetilde{\Talg}'^F}(\theta')\in\Irr
\left(\widetilde{\Talg}'^{\circ
F}\right)^{\tau}$ and $\theta_0''=\Res_{\widetilde{\Talg}''^{\circ
F}}^{\widetilde{\Talg}''^F}(\theta'')\in\Irr\left(\widetilde{\Talg}''^{\circ
F}\right)^{\tau}$. We have
$$\cyc{R_{\widetilde{\Talg}'}^{\widetilde{\Galg}}(\theta),R_{\widetilde{\Talg}''}^{\widetilde{\Galg}}(\theta'')}_{\tau}=0,$$ if~$(\widetilde{\Talg}'^{\circ},\theta_0')$ and~$(\widetilde{\Talg}''^{\circ},\theta_0'')$
are not~$\left((\widetilde{\Galg}^{\tau})^{\circ}\right)^F$-conjugate and
$$\cyc{R_{\widetilde{\Talg}'}^{\widetilde{\Galg}}(\theta'),R_{\Talg'}^{\widetilde{\Galg}}(\theta')}_{\tau}= \frac{1}{|\left(\widetilde{\Talg}'^{\tau}\right)^{\circ F}|}
\left|\{n\in\Nn_{(\widetilde{\Galg}^{\tau})^{\circ F}}(\widetilde{\Talg}')\ |\ ^n\theta_0'=\theta_0'\}\right|.$$
Here, $\cyc{\psi_1,\psi_2}_{\tau}$
denotes~$\frac{1}{|\Galg^F|}\sum\limits_{g\in\Galg^{F}}\psi_1(g\tau)\overline{\psi_2(g\tau)}$.
\end{proposition}

\subsection{Generalized Springer correspondence}

The generalized Springer correspondence is extended to disconnected
groups in~\cite[II]{Lustdisconnected}. It is extended for the
connected component~$D=\Galg\tau$ of~$\widetilde{\Galg}$:
since~$\tau$ is unipotent, we can associate to every pair~$(u,\phi)$
of~$\mathcal{N}_{D}$ ($u$ is a unipotent element lying in~$D$
and~$\phi$ is an irreducible character of~$A(u)$ defined in~\S\ref{secconjclass})
a~$4$-tuple~$(\Lalg,v,\psi,\rho)$ as in the connected case,
where~$\Lalg$ is a~$\tau$-stable Levi subgroup of~$\Galg$,~$(v,\psi)$ is a cuspidal 
pair of~$\mathcal{N}_{\Lalg\tau}$
and~$\rho$ is an irreducible character of~$\Nn_{\Galg^{\tau}}(\Lalg)/\Lalg$;
see~\cite[II.8.8,II.11.10]{Lustdisconnected}. For the pairs corresponding to~$\Talg$ 
(denoted~$\mathcal{N}_{D,0}$ as previously), the generalization to the disconnected case
was independently introduced by Sorlin~\cite{Sorlin}.

\begin{proposition}\label{spdisconnected} Let~$\Galg_0$, $\Galg$
and~$\widetilde{\Galg}$ as in~\S\ref{rescal}. Let~$u$ be a
unipotent element of~$\Galg_0$ and~$(u,\phi)\in\mathcal{N}_0$ and 
let~$\rho_{u,\phi}\in W_0$ be the Springer correspondent of~$(u,\phi)$.
Then~$(f(u),\phi\circ\mu^{-1})\in\mathcal{N}_{D,0}$ (where~$f$ and~$\mu$ are the maps defined in Proposition~\ref{bijection}
) and the Springer correspondent
corresponding to it is~$\rho_{u,\phi}$ (viewed as a character
of~$\Nn_{\Galg^{\tau}}(\Talg)/\Talg\simeq W_0$).
\end{proposition} 

\begin{preuve}{}
This is a consequence of Proposition~\ref{bijection} and~\cite[Theorem 3.3]{MalleSorlin}, applied with~$\sigma$ being trivial.

\end{preuve}

\begin{corollaire}
The Springer correspondence for the pairs lying in~$\mathcal{N}_{D,0}$ 
for the disconnected groups 
of type~$B_2$ and~$F_4$ are given in Tables~\ref{SPB2} and~\ref{SPF4}
respectively. In Table~\ref{SPB2} (resp.~Table~\ref{SPF4}), 
the representatives~$u_i$ (resp.~$x_i$) are now the images under~$f$ defined in 
Proposition~\ref{bijection} of the corresponding elements of~$\Galg_0$.
\end{corollaire}

\begin{remarque} In fact, for the disconnected cases of type~$B_2$
and~$F_4$, the generalized Springer correspondence of the pairs lying
in~$\mathcal{N}_D$ is in correspondence with the generalized
correspondence of the pairs lying in~$\mathcal{N}$.
The generalized Springer correspondence for disconnected groups
indeed satisfies the same properties as the generalized Springer
correspondence for connected groups. These properties are recalled
in~\cite[\S2]{Mallesp}. Moreover in~\cite{spSP}, Spaltenstein proves that
the generalized Springer correspondence for a simple group of type~$F_4$
in characteristic~$2$ corresponding to the pairs which are not in~$\mathcal{N}_0$
can be computed using only these properties.
\end{remarque}

\subsection{The disconnected case of type~$B_2$}

We use the notation of~\cite{brunat}. Note that the symbol~$q$ in~\cite{brunat}
is the symbol~$q^2$ in this paper. The group~$\Galg^F$ has~$4$
unipotent characters~$1_{\Galg^F}$,~$\theta_1$,~$\theta_4$
and~$\theta_5$ which extend to~$\widetilde{\Galg}^F$. If~$\theta$ is
such a character, we denote by~$\widetilde{\theta}$ the extension
such that~$\widetilde{\theta}(\tau)>0$. The values of these
extensions are given in~\cite[Table 9]{brunat}. Since every class
of~$\widetilde{\Galg}$ lying in the coset~$\Galg\tau$ is real, it
follows from Remark~\ref{shin} that the class of~$f(x)$ (where $x\in\Galg_0$ 
and~$f$ defined in Proposition~\ref{bijection}) is the Shintani
correspondent of~$x$ in~$\Galg_0^{F_0^2}\semi{F_0}$ when we use the
identifications of~Proposition~\ref{isogroupe}. This correspondence is
explicitly described in~\cite[Proposition 4.1]{brunat}.
Following Proposition~\ref{tores}, the $\widetilde{\Galg}^F$-classes of 
rational "tori" of~$\widetilde{\Galg}$ 
are parametrized by the~$F_0$-classes of~$W=W(\Galg_0)$ and
a system of representatives is given in Table~\ref{WB2}. For
such a representative~$w$, we denote by~$R_w$ the restriction
of~$R_{\widetilde{\Talg}_w}^{\widetilde{\Galg}}(1_{\widetilde{\Talg}_w^F})$ 
to~$\Galg^F\tau$.

\begin{proposition}
We have on the coset~$\Galg^F\tau$:
$$\begin{array}{lcl}
R_1&=&1_{\widetilde{\Galg}}+\widetilde{\theta}_4\\
R_{w_a}&=&1_{\widetilde{\Galg}}-\widetilde{\theta}_1-\widetilde{\theta}_5-\widetilde{\theta}_4\\
R_{w_aw_bw_a}&=&1_{\widetilde{\Galg}}+\widetilde{\theta}_1+\widetilde{\theta}_5-\widetilde{\theta}_4\
\end{array}$$

\end{proposition}

\begin{preuve}{}
Using the proof of~\cite[Proposition 3.1]{brunat}, we deduce
that~$R_1=1_{\widetilde{\Galg}}+\widetilde{\theta}_4$ on~$\Galg^F\tau$.
Let~$w\in W$; note that~$\cyc{R_w,1_{\widetilde{\Galg}}}_{\tau}=1$, where~$\cyc{,}_{\tau}$ is the restricted scalar product introduced in Proposition~\ref{propDM}.
We indeed have~$$\cyc{R_{\widetilde{\Talg}_w}^{\widetilde{\Galg}}(1_{\widetilde{\Talg}_w^F}),1_{\widetilde{\Galg}}}_{\tau}=\cyc{1_{\widetilde{\Talg}_w^F},^*\!\!R_{\widetilde{\Talg}_w}^{\widetilde{\Galg}}(1_{\widetilde{\Talg}_w^F})}_{\Talg_w^{\circ F}\tau}=\cyc{1_{\widetilde{\Talg}_w^F},1_{\widetilde{\Talg}_w^F}}_{\Talg_w^{\circ F}\tau}=1,$$
where~$^*\!R_{\widetilde{\Talg}_w}^{\widetilde{\Galg}}$ is the adjoint of~$R_{\widetilde{\Talg}_w}^{\widetilde{\Galg}}$.
Using Proposition~\ref{propDM}, we deduce that
$$\begin{array}{lcl}\cyc{R_{w_a},R_{w_a}}_{\tau}&=&\cyc{R_{w_aw_bw_a},R_{w_aw_bw_a}}_{\tau}=4,\\
\cyc{R_1,R_{w_a}}_{\tau}&=&\cyc{R_1,R_{w_aw_bw_a}}_{\tau}=0.
\end{array}
$$
We then deduce that
$$\cyc{\widetilde{\theta}_4,R_{w_a}}_\tau=\cyc{\widetilde{\theta}_4,R_{w
_aw_bw_a}}_\tau=-\cyc{1_{\widetilde{\Galg}},R_{w_a}}_\tau=-1.$$
Moreover, since the constituents of~$R_{w_a}$ and~$R_{w_aw_bw_a}$
are extensions of unipotent characters of~$\Galg^F$, we
deduce that the other constituents of these characters
are~$\widetilde{\theta}_1$ and~$\widetilde{\theta}_5$ (possibly
with multiplicity~$0$). However, the~$\tau$-norm of these
characters is~$4$. It then follows that~$\widetilde{\theta}_1$
and~$\widetilde{\theta}_5$ occur with multiplicity~$\pm
1$. Now, using~\cite[Theorem 4.13]{DMnonconnexe}, we have
$$R_{w_a}(\tau)=R_{\Talg_{0,w_a}}^{\Galg_0}(1)=-(q^2-1)(q^2+\sqrt{2}q+1)
,$$ and similarly~$R_{w_aw_bw_a}(\tau)=-(q^2-1)(q^2-\sqrt{2}q+1)$.\\
Since~$\widetilde{\theta}_1(\tau)=\widetilde{\theta}_5(\tau)=\frac{1}
{\sqrt{2}}q(q^2-1)$, the result holds.

\end{preuve}
\begin{corollaire}\label{cor}
We choose the extensions of the~$F$-stable characters
of~$W^{\tau}$ to~$W^{\tau}\semi{F}$ as in Table~\ref{WB2}. We then have
$$
\begin{array}{lcl}
R_{\widetilde{1}}&=&1_{\widetilde{\Galg}^F},\\
R_{\widetilde{\varepsilon}}&=&\widetilde{\theta}_4,\\
R_{\widetilde{\chi}}&=&\frac{\sqrt{2}}{2}(\widetilde{\theta}_1+\widetilde{\theta}_5).
\end{array}
$$
Moreover, the values of these class functions on the unipotent elements
lying in the coset~$\Galg^F\tau$ are given in Table~\ref{almostB2},
where we replace $u_1$, $u_2$,~$\rho$ and~$\rho^{-1}$ by the
elements denoted by~$(1,\s)$, $(x_{a+b},\s)$, $(x_a,\s)$
and~$(x_ax_{a+b},\s)$ in~\cite{brunat} respectively. 
\end{corollaire}

Let~$\mathcal{N}_{D,0}^F$ be the set of pairs of~$\mathcal{N}_{D,0}$ which are $F$-stable. When~$(u,\phi)\in\mathcal{N}_{D,0}^F$,
we can associate to~$(u,\phi)$ a class function~$Y_{u,\phi}$ on the unipotent elements
of~$\Galg^F\tau$; see~\cite[IV.19.6, IV.19.8]{Lustdisconnected}. When~$u$ is a unipotent element of~$\Galg^F\tau$, 
we denote by~$\mathcal{B}_u$ the variety of Borel subgroups of~$\Galg$ which are~$u$-stable (for the conjugation in~$\widetilde{\Galg}$).
We now have:

\begin{itemize}
\item A Springer correspondence for the pairs lying in~$\mathcal{N}_{D,0}^F$;
\item A basis of the space of class functions defined on the unipotent elements of~$\Galg^F\tau$.
\item The integers~$d_u$ giving the dimensions of the varieties~$\mathcal{B}_u$.
\end{itemize} 
We can then formally apply  Lusztig's
Algorithm with these data as input. Remark that if~$u$ is a unipotent
element of~$\Galg_0$, then the varieties~$\mathcal{B}_u$ and~$\mathcal{B}_{f(u)}$ (with~$f$ defined
in Proposition~\ref{bijection}) are isomorphic. The isomorphism is given by
$$\Balg'\in\mathcal{B}_u\mapsto(\Balg',\Balg')\in\mathcal{B}_{f(u)}.$$
We have seen in Proposition~\ref{bijection} that~$f(u)$
lies in a rational unipotent class of~$\Galg\tau$. Let~$u'$ be
a rational element in the class of~$f(u)$ obtained by the process
of the proof of Proposition~\ref{bijection}. We have shown that the action
of~$F$ on~$A(u')$ is the same as
the one of~$F$ on~$A(f(u))$. We thus obtain a bijection
between the pairs~$(u',\phi')\in\mathcal{N}_{D,0}^{F_0}$ and the
pairs~$(f(u),\phi\circ\mu^{-1})$, with~$\phi\circ\mu^{-1}$ being an
$F$-stable character of~$A(f(u))$. Since the action
of~$F$ on~$A(f(u))$ is the same as the one of~$F_0$
on~$A_{\Galg_0}(u)$, it follows that~$\mathcal{N}_0^{F_0}$
and~$\mathcal{N}_{D,0}^F$ are in bijection. 
%
Using the proof of~\cite[IV,~Lemma~9.7]{Lustdisconnected},
we define the class functions~$\varphi_{f(u),\phi}$
for~$(u,\phi)\in\mathcal{N}_{0}^{F_0}$ by the class function on
the set of unipotent elements of~$\Galg^F\tau$ as follows:
the conjugacy classes of~$\widetilde{\Galg}^F$ lying in~$\Cl(f(u)^F)$
are in bijection with the~$F$-classes
of~$A(f(u))$ (see the proof of Proposition~\ref{bijection}).
If~$a$ is a representative of such a class, we denote by~$\widetilde{a}$ a
representative of the corresponding class of~$\widetilde{\Galg}^F$.
Let~$g\in\widetilde{\Galg}^F$, we set
$$\varphi_{f(u),\phi}(g)=\left\{\begin{array}{ll}
\widetilde{\phi}(aF)&\textrm{if }g\ \textrm{is
conjugate in } \widetilde{\Galg}^F\ \textrm{to } \widetilde{a},\\
0&\textrm{otherwise.}
\end{array}\right.$$
%
Here~$\widetilde{\phi}$
denotes an extension of~$\phi\circ\mu^{-1}$ to~$A(f(u))\semi{F}$.
We choose for~$\widetilde{\phi}$ the same extension as in
Table~\ref{YB2}. For the notation, we replace in Table~\ref{YB2} the
functions~$\varphi_{u,\phi}$ by~$\varphi_{f(u),\phi}$ and the classes as in
Corollary~\ref{cor}.

\begin{proposition}\label{propB2}
The outputs obtained using Lusztig's algorithm with these data
are precisely the uniform almost characters of~$\Galg^F\tau$.
\end{proposition}

\subsection{Conjecture for the disconnected case of type~$F_4$}

Since Lusztig's algorithm is valid for the disconnected case of type~$B_2$, 
we can suppose that this will also be the case for the disconnected type~$F_4$. 
We state a conjecture of the same type as Proposition~\ref{propB2}.

\begin{conjecture} Suppose~$\Galg_0$ is a simple group of type~$F_4$ and~$F_0$ 
the generalized Frobenius map that defines a Ree group of type~$F_4$. 
Then the values of the
uniform almost characters of the group~$\Galg_0^{F_0^2}\semi{F_0}$ 
on the coset~$\Galg_0^{F_0^2}F_0$ are given in Table~\ref{almostF4}, 
when we replace~$u_i$ by its Shintani correspondent~$N_{F_0/F_0^2}(u_i)$.
\end{conjecture}

\noindent{\bf Acknowledgment}\\

\noindent I wish to thank heartily C\'edric Bonnaf\'e to orient me in this
direction, and I wish to thank deeply Frank L\"ubeck and Gunter Malle
for the many valuable discussions that they granted to me.


\begin{thebibliography}{0}
\bibitem{bonnafe}
{\sc C. Bonnaf\'e}, Sur les caract\`eres des groupes r\'eductifs finis \`a centre non connexe: applications aux groupes sp\'eciaux lin\'eaires et unitaires, to appear in Ast\'erisque.
\bibitem{brunat}
{\sc O. Brunat}, The Shintani descents of Suzuki groups and their consequences, Journal of Algebra 303, 2006, 869--890.
\bibitem{brunat1}
{\sc O. Brunat}, On the extension of~$\Gf_2(3^{2n+1})$ by the exceptional graph automorphism, submitted for publication.
\bibitem{Carter}
{\sc R.W. Carter}, Simple groups of Lie Type, Wiley, New York,
1972.
\bibitem{Carter2}
{\sc R.W. Carter}, Finite Groups of Lie Type: Conjugacy classes and
Complex Characters, Wiley, New York, 1985.
\bibitem{Digne}
{\sc F. Digne}, Descente de Shintani et restriction des scalaires, Journal of London Math. Society, vol. 59, No 3 (1999) 
\bibitem{DM}
{\sc F. Digne and J. Michel}, Fonctions $\mathcal{L}$ des
vari\'et\'es de Deligne-Lusztig et descente de Shintani, Bull. S.M.F., m\'emoires {\bf 20}, 1985.
\bibitem{DMnonconnexe}
{\sc F. Digne and J. Michel}, Groupes r\'eductifs non connexes, Ann. \'Ecole Normale Sup\'erieure, 4i\`eme s\'erie, t. 27, 1994, 345--406.
\bibitem{gap}
{\sc The Gap Group}, Groups, Algorithms, and Programming,
      Version 4.4.7; http://www.gap-system.org, 2006.
\bibitem{GM}
{\sc M. Geck and G. Malle}, Fourier transforms and Frobenius
eigenvalues for finite Coxeter groups, Journal of Algebra 260, 2003.
162--193.
\bibitem{Kondo}
{\sc T. Kondo}, The characters of the Weyl group of type~$F_4$, J. Fac. Sci. Univ. Tokyo, 11, 1965, 145--153.
\bibitem{Lust}
{\sc G. Lusztig}, {Characters of reductive groups over a finite
field}, Annals Math.\ Studies 107, Princeton University Press,
1984.
\bibitem{Lustinter}
{\sc G. Lusztig}, Intersection cohomology complexes on reductive group, Invent.
Math., 75, 1984, 205--272.
\bibitem{Lustsheave}
{\sc G. Lusztig}, Character sheaves, Adv. Math. 56 (1985), 193-237; 
II, Adv. Math. 57 (1985) 226-265; III, Adv. Math. 57 (1985), 266-315; IV, Adv. Math. 59 (1986), 1-63; V, Adv. Math. 61 (1986), 103--155.
\bibitem{Lustdisconnected}
{\sc G. Lusztig}, Character sheaves on disconnected groups II, Represent. Theory 8 (2004), 72--124; IV, Represent. Theory 8 (2004), 145--178.
\bibitem{Malle2F4}
{\sc G. Malle}, Die unipotenten Charaktere von $^2F_4(q^2)$, Comm. Algebra, 18, 1990, 2361--2381.
\bibitem{Malle1}
{\sc G. Malle}, Generalized Deligne-Lusztig Characters, Journal of Algebra 159, 1993, 64--97.
\bibitem{Mallesp}
{\sc G. Malle}, Springer correspondence for disconnected exceptional groups, Bull. London Math. Soc. 37, 2005, 391--398.
\bibitem{MalleSorlin}
{\sc G. Malle and K. Sorlin}, Springer correspondence for disconnected groups, Math. Z. 246, 2004, 291--319.
\bibitem{Shinoda2}
{\sc K. Shinoda}, The conjugacy classes of Chevalley groups of type
$(F_4)$ over finite fields of characteristic~$2$, J. Fac. Sci. Univ.
Tokyo, 21, 1974, 133-159. 
\bibitem{Shinoda}
{\sc K. Shinoda}, The conjugacy classes of the finite Ree groups of type $(F_4)$, J. Fac. Sci. Univ. Tokyo, 22, 1975, 1--15.
\bibitem{shojiSC}
{\sc T. Shoji}, Geometry of orbits and Springer correspondence, Soc. Math. de France, Ast\'erisque~168, 1988, 61--140.
\bibitem{shoji}
{\sc T. Shoji}, Character sheaves and almost characters of reductive groups, Adv Math. 111 (1995), 244--313;
\bibitem{Sorlin}
{\sc K. Sorlin}, Springer correspondence in non connected reductive groups, J. Reine Angew. math. 508, 2004, 197--234.
\bibitem{Sp}
{\sc N. Spaltenstein}, Classes unipotentes et sous-groupes de Borel, Lecture Notes in Math. 
946, Springer, Berlin Heidelberg New York, 1982
\bibitem{spSP}
{\sc N. Spaltenstein}, On the generalized Springer correspondence for exceptional groups, Algebraic groups and related topics, Adv. Stud. Pure Math. 6 (kinokuniya and North-Holland, Tokyo and Amsterdam, 1985), 317--338.
\bibitem{Suzuki}
{\sc M. Suzuki}, On a class of doubly transitive groups. I, Annals
of Math. 75, 1962, 105--145.
\bibitem{Ward}
{\sc H.N. Ward}, On Ree's series of simple groups, Trans. Amer.
Math. Soc., 121, 1966, 62--89.
\end{thebibliography}
\end{document}